\documentclass[11pt,twoside]{article}

\usepackage[cm]{styleset}

\usepackage{macros}

\title	 {\normalfont\headfam\itshape\LARGE
          Four-manifold invariants from higher-rank bundles}
\author	 {P. B. Kronheimer}
\address {Harvard University, Cambridge MA 02138}

\begin{document}

\pagestyle{numbered}	
\maketitle

\section{Introduction}

\subsection{Background}

Donaldson's polynomial invariants for a smooth $4$-manifold $X$ are
defined in \cite{Donaldson-polynomials} using the moduli spaces of
anti-self-dual connections on a vector bundle $E \to X$. The structure
group of the bundle $E$ is generally taken to be either $\SU(2)$ or
$\SO(3)$, both in \cite{Donaldson-polynomials} and in the subsequent
mathematical literature, even though much of the underlying analysis
is applicable to bundles with arbitrary compact structure group.  In
the physics literature, there is the work of Mari\~no and Moore, who
argue in \cite{Marino-Moore} that the $\SU(2)$ Donaldson invariants
have a well-defined generalization for $\SU(N)$ bundles. They also
propose a formula for the $\SU(N)$ invariants when $X$ has simple type, which
suggests that these higher-rank invariants contain no new topological information.
The purpose of this paper is to give a mathematical definition of the
$\SU(N)$ Donaldson invariants, and to calculate some of them in an
interesting family of examples arising from knot complements. The
results we obtain are consistent with the predictions of
\cite{Marino-Moore}.

We recall some background material on the Donaldson invariants.  The
notion of ``simple type'' was introduce in \cite{KM-Structure}.  For
$4$-manifolds $X$ with this property, the Donaldson invariants arising
from the Lie groups $\SU(2)$ and $\SO(3)$ were shown to be determined
by a rather small amount of data.  To state this result, we recall
that, in the case that $E$ is an $\SU(2)$ bundle with $c_{2}(E)[X] =
k$, the corresponding polynomial invariant $q_{k}$ is defined
whenever $b_{2}^{+}(X) - b_{1}(X)$ is odd and $b_{2}^{+}(X)$ is at
least $2$; the invariant is a homogeneous polynomial function
\[ 
     q_{k} : H_{2}(X;\Z) \to \Z
\]
of degree $d/2$, where $d = 8k - 3(b^{+}_{2} - b_{1} + 1)$ is the
expected dimension of the instanton moduli space of instanton number
$k$ on $X$. If $X$ has simple type and $b_{1}(X) = 0$, then the structure
theorem from \cite{KM-Structure} establishes that there is a
distinguished collection of $2$-dimensional cohomology classes $K_{i}
\in H^{2}(X;\Z)$ ($i=1,\dots, r$) and rational coefficients $a_{i}$
(independent of $k$) such that $q_{k}/(d/2)!$ is the term of degree
$d/2$ in
the Taylor expansion of the function
\[
\begin{gathered}
    \mathcal{D}_{X} : H_{2}(X;\R) \to \R  
\end{gathered}
\]
given by
\begin{equation}\label{eq:Dseries}
            \mathcal{D}_{X}(h) = \exp \left(\frac{Q(h)}{2}\right) \sum_{i=1}^{r} a_{i}
            \exp {K_{i}(h)}.
\end{equation}
Here $Q$ is the intersection form on $H_{2}(X;\Z)$.
There is a
similar formula for the $\SO(3)$ invariants, involving the same
``basic classes'' $K_{i}$ and coefficients $a_{i}$, which means
that there is no more information in the $\SO(3)$ invariants than is
already contained in the $\SU(2)$ invariants.

Witten's conjecture \cite{Witten} states that for simply-connected 4-manifolds of
simple type, the coefficients $a_{i}$ and basic classes $K_{i}$ are
determined by the Seiberg-Witten invariants of $X$: we should have
\[
            a_{i} = 2^{2+(7\chi(X) + 11\sigma(X))/4}\mathit{SW}_{X}(\mathfrak{s}_{i})
\]
and
\[
           K_{i} = c_{1}(S^{+}_{\mathfrak{s}_{i}}),
\]
where the $\mathfrak{s}_{i}$ are the spin-c structures on $X$ for which
the Seiberg-Witten invariant $\mathit{SW}_{X}(\mathfrak{s}_{i})$ is
non-zero  and $S^{+}_{\mathfrak{s}_{i}} \to X$ is the half-spin bundle
corresponding to $\mathfrak{s}_{i}$. This conjecture was extended to
the case of $\SU(N)$ Donaldson invariants in \cite{Marino-Moore}:
the conjectured formula for the $\SU(N)$ invariants
is a homogeneous expression of degree $N-1$ in the quantities
$\mathit{SW}_{X}(\mathfrak{s}_{i})$.

\subsection{Summary of results}

In this paper 
we shall give a quite careful construction of the $\SU(N)$ and
$\PSU(N)$ Donaldson
invariants for $4$-manifolds $X$ with $b^{+}_{2} \ge 2$. As is the case
already with $N=2$, the $\PSU(N)$ invariants are best defined by
choosing a bundle $P \to X$ with structure group $U(N)$ (and
possibly non-zero first Chern class) and then studying a moduli
space of connections $A$ in the associated adjoint bundle.
(Note that not every principal $\PSU(N)$ bundle can be lifted to a
$U(N)$ bundle, so there is a slight loss of generality at this point.)

The new difficulties that arise in the higher-rank case are two-fold.
First, we can expect singularities in the moduli space corresponding
to \emph{reducible} anti-self-dual connections. We shall bypass this
difficulty initially by considering only the case that $c_{1}(P)$ is
not divisible by any prime dividing $N$ in the lattice
$H^{2}(X;\Z)/\mathrm{torsion}$. In this case, if $b^{+}_{2}(X)$ is
non-zero, then for a generic Riemannian metric on $X$ there will be no
reducible solutions to the equations, as we explain in
section~\ref{subsec:no-reducibles}.  If $c_{1}(P)$ is not coprime to
$N$ in this way, we can replace $X$ with the blow-up $\tilde X = X \#
\bar{\CP}^{2}$, and then define the invariants of $X$ using moduli
spaces on $\tilde X$. (This is the approach taken in
\cite{KM-Structure} for the $\SU(2)$ invariants.)  The second
difficulty is that for $N>2$, one does not know that the irreducible
parts of the $\SU(N)$ or $\PSU(N)$ moduli spaces will be smooth for
generic choice of Riemannian metrics: we do not have the ``generic
metrics theorem'' of Freed and Uhlenbeck \cite{Freed-Uhlenbeck}, and
the anti-self-duality equations must be artificially perturbed to
achieve transversality.

Having defined the higher-rank invariants,
we shall also calculate them in the case corresponding to a
zero-dimensional moduli space on a particular family of
$4$-manifolds.  Let $K$ be a knot in $S^{3}$, let $M$ be the knot
complement (a compact $3$-manifold with torus boundary), let $X$ be a
$K3$ surface with an elliptic fibration, and let $T\subset X$ be a
torus arising as one of the fibers.  The following construction is
considered by Fintushel and Stern in \cite{Fintushel-Stern-Knot}. Let $X^{o}$ denote
the $4$-manifold with boundary $T^{3}$ obtained by removing an open
tubular neighborhood of $T$, and let $X_{K}$ be the closed
$4$-manifold formed as
\[
            X_{K} = (S^{1}\times M) \cup_{\phi}  X^{o},
\]
where $\phi : S^{1} \times \partial M \to \partial X^{o}$ is an
identification of the two $3$-torus boundaries, chosen so that $\phi$
maps $(\mathrm{point})\times (\mathrm{longitude})$ to a curve on $\partial
X^{o} \subset X$ which is the boundary of a $2$-disk transverse to
$T$. It is shown in \cite{Fintushel-Stern-Knot} that $X_{K}$ has the same homotopy
type as $X$, but that the Seiberg-Witten invariants of $X_{K}$
encode the Alexander polynomial of $K$. If Witten's conjecture holds,
then the function $\mathcal{D}$ that determines the $\SU(2)$ Donaldson
invariants (see \eqref{eq:Dseries}) is given by
\[
        \mathcal{D}_{X_{K}} =
        \exp \left(\frac{Q(h)}{2}\right) \Delta_{K}( \exp( 2 F(h))),
\]
where $F$ is the cohomology class Poincar\'e dual to the torus fibers
and
\[
            \Delta_{K}(t) = \sum_{k=-n}^{n} a_{k} t^{k} 
\]
is the symmetrized Alexander polynomial of the knot $K$.

We shall consider a particular $U(N)$ bundle $P\to X_{K}$
with $(c_{1}(P)\cdot F) = 1$ and $c_{2}(P)$ chosen so that the
corresponding moduli space is zero-dimensional. The zero-dimensional
moduli space gives us an integer-valued invariant of $X_{K}$ (the
signed count its points), which we are able to calculate when $N$ is
odd. The answer we obtain is the integer
\[
       \prod_{k=1}^{N-1} \Delta_{K}({e^{2\pi i k/N}}),
\]
as long as this quantity is non-zero. It seems likely that this answer
is correct also when $N$ is even, up to an overall sign (the answer is consistent with
Witten's conjecture when $N=2$ for example); but we have not overcome
an additional difficulty that arises in the even case: see
Lemma~\ref{lem:same-sign} and the remarks at the end of
section~\ref{subsec:KnotComplements}.

\medskip

\subparagraph{Acknowledgement.}

This paper owes a lot to the author's continuing collaboration with Tom
Mrowka. In particular, the discussion of holonomy perturbations in
section~\ref{sec:perturb} is based in part on earlier, unpublished
joint work.

\section{Gauge theory for $\SU(N)$}

\subsection{The configuration space}

Let $X$ be a smooth, compact, connected, oriented $4$-manifold. We
consider
the Lie group $U(N)$ and its adjoint group, $\PU(N)$ or
$\PSU(N)$.
Let $P \to X$
be a smooth principal $U(N)$-bundle, and let $\adP$ denote the associated
principal bundle with structure group the adjoint group.  Fix an integer $l \ge
3$, and let $\cA$ denote the space of all connections in
$\adP$  of
Sobolev class $L^{2}_{l}$. The assumption $l\ge 3$ ensures that
$L^{2}_{l}\subset C^{0}$: we will use this assumption in
Lemma~\ref{lem:transverse}, although at most other points $l\ge 2$
would suffice.
If we fix a particular connection $A_{0}$,
we can describe $\cA$ as
\[
        \cA = \{\, A_{0} + a \mid a \in L^{2}_{l}(X; \g_{P} \otimes
        \Lambda^{1}X) \,\},
\]
where $\g_{P}$ is the vector bundle with fiber $\su(N)$ associated to
the adjoint representation.  The \emph{gauge group} $\G$ will be the
group of all automorphisms of $P$ of Sobolev class $L^{2}_{l+1}$ which
have determinant 1.  That is, $\G$ consists of all
$L^{2}_{l+1}$ sections of the fiber bundle $\SU(P) \to X$ with fiber
$\SU(N)$, associated to the adjoint action of $\PSU(N)$ on $\SU(N)$. 
The group $\G$ acts on $\cA$, but the action is not effective. The
subgroup that acts trivially is a cyclic group $Z \subset \G$ that we
can identify with the center of $\SU(N)$.

\begin{definition}
    For a connection $A \in \cA$, we write $\Gamma_{A}$ for the
    stabilizer of $A$ in $\G$. We say that $A$ is \emph{irreducible}
    if $\Gamma_{A} = Z$.
\end{definition}

We write $\cA^{*}\subset \cA$ for the irreducible connections. We
write $\B$ and $\B^{*}$ for the quotients of $\cA$ and $\cA^{*}$ by
the action of  $\G$. The space $\B^{*}$ has the structure of a Banach
manifold in the usual way, with coordinate charts being provided by
slices to the orbits of the action of $\G$. We write $[A] \in \B$ for
the gauge-equivalence class represented by a connection $A$.

A
$U(N)$-bundle $P$ on $X$ is determined up to isomorphism by its first and second
Chern classes. We shall write $w \to X$ for the line bundle associated
to $P$ by the $1$-dimensional representation $g \mapsto \det(g)$ of
$U(N)$, so that $c_{1}(w) = c_{1}(P)$. We set
\begin{equation}\label{eq:kappa-def}
\begin{aligned}
\k &= \left\langle c_{2}(P) - ((N-1)/(2N)) c_{1}(P)^{2} , [X]
   \right\rangle \\
   &=  -(1/2N) \left\langle p_{1}(\g_{P}) , [X]
   \right\rangle .
   \end{aligned}
\end{equation}
The above expression is constructed so
that $\kappa$ depends only on $\adP$, and so that, in the case that
$X$ is $S^{4}$, the possible values of $\kappa$ run through all
integers as $P$ runs through all isomorphism classes of $U(N)$-bundles.
We may label the configuration spaces $\B$ etc.~using $w$ and $\k$
to specify the isomorphism class of the bundle. So we have
$\B^{w,*}_{\k} \subset \B^{w}_{\k}$ and so on. We refer to $\kappa$ as
the \emph{instanton number}, even though it is not an integer.

There is a slightly different viewpoint that one can take to describe
$\cA$. Fix a connection $\theta$ in the line bundle $w$. Then for each
$A$ in $\cA$, there is a unique $U(N)$-connection $\tilde A$ in the
$U(N)$-bundle $P$ such that the connections which $A$ induces in $\adP$
and $w$ are $A$ and $\theta$ respectively. Thus we can identify $\cA$
with the space $\tilde\cA$ of $U(N)$-connections inducing the connection
$\theta$ in $w$.

This viewpoint is helpful in understanding the
stabilizer $\Gamma_{A}$. Pick a basepoint $x_{0} \in X$, and let $A$
be a connection in $\adP$. An element $g$ in $\Gamma_{A}$ is
completely determined by its restriction to the fiber over $x_{0}$,
so we can regard $\Gamma_{A}$ as a subgroup of the copy of a copy of
$\SU(N)$, namely $\SU(P)_{x_{0}}$. Let $\tilde A$ be the
$U(N)$-connection in $P \to X$, corresponding to $A$ and $\theta$ as
above, and let $\tilde H$ be the group of automorphisms of $P_{x_{0}}$
generated by the holonomy of $\tilde A$ around all loops based at
$x_{0}$. Then $\Gamma_{A}$ is the centralizer of $\tilde H$ in
$\SU(P)_{x_{0}}$.  The structure of $\Gamma_{A}$ can then be described
as follows. For any isomorphism of inner product spaces
\[
         \psi:    \bigoplus_{i=1}^{r} \left( \C^{n_{i}} \otimes\C^{m_{i}}
            \right) \to \C^{N},
\]
we get a representation
\[
        \rho : \prod_{i=1}^{r} U(n_{i}) \to U(N).
\]
The group $\Gamma_{A} \subset \SU(P)_{x_{0}} \cong \SU(N)$ is
isomorphic to
\[
    \rho\left( \prod_{i=1}^{r} U(n_{i}) \right) \cap \SU(N)
\]
for some such $\psi$. In particular, we make the important observation
that $\Gamma_{A}$ has positive dimension unless $r=1$ and $n_{1}=1$:
that is to say, $\Gamma_{A}$ has positive dimension unless $A$ is
irreducible.

\subsection{The moduli space}

Let $M^{w}_{\kappa} \subset \B^{w}_{\kappa}$ be the space of
gauge-equivalences classes of anti-self-dual connections:
\[
   M^{w}_{\kappa} = \{ \, [A] \in B^{w}_{\k} \mid F^{+}_{A} = 0 \, \}.
\]
If $A$ is anti-self-dual, then it is gauge-equivalent to a smooth
connection. Supposing that $A$ itself is smooth, we have a complex
\[
        \Omega^{0}(X;\g_{P}) \stackrel{d_{A}}{\to}
        \Omega^{1}(X;\g_{P}) \stackrel{d^{+}_{A}}{\to}
        \Omega^{+}(X;\g_{P}) .
\]
We denote by $H^{0}_{A}$, $H^{1}_{A}$ and $H^{2}_{A}$ the cohomology
groups of this elliptic complex. If $A$ is not smooth but only of
class $L^{2}_{l}$, we can define the groups $H^{i}_{A}$ in the same
way after replacing the spaces of smooth forms above by Sobolev
completions.

\begin{lemma}
    The vector space $H^{0}_{A}$ is zero if $A$ and only if $A$ is irreducible. 
\end{lemma}

\begin{proof}
    In general, $H^{0}_{A}$ is the Lie algebra of $\Gamma_{A}$. We
    have already observed that $\Gamma_{A}$ has positive dimension if
    $A$ is reducible.
\end{proof}

We say that $A$ is
\emph{regular} if $H^{2}_{A}$ vanishes.  If $A$ is regular and
irreducible, then
$M^{w}_{\kappa}$ is a smooth manifold in a neighborhood of $[A]$, with
tangent space isomorphic to $H^{1}_{A}$. The dimension of the moduli
space near this point is then given by minus the index of the complex
above, which is
\begin{equation}\label{eq:index}
    d = 4N \kappa - (N^{2}-1)(b^{+}_{2}(X) - b_{1}(X) + 1).
\end{equation}
We say that the moduli space is regular if $A$ is regular for all
$[A]$ in $M^{w}_{\k}$.

Uhlenbeck's theorem allows us to compactify the moduli space
$M^{w}_{\kappa}$ in the usual way. We define an ideal anti-self-dual
connection of instanton number $\kappa$ to be a pair $([A], \bx)$,
where $A$ is an anti-self-dual connection belonging to
$M^{2}_{\kappa-m}$, and $\bx$ is a point in the symmetric product
$X^{m}/S_{m}$. The space of ideal anti-self-dual connections is
compactified in the usual way, and we define the Uhlenbeck
compactification of $M^{w}_{\kappa}$ to be the closure of the
$M^{w}_{\kappa}$ in the space of ideal connections:
\[
            \bar{M}^{w}_{\kappa} \subset \bigcup_{m} M^{w}_{\kappa-m}
            \times (X^{m}/S_{m}).
\]

\subsection{Avoiding reducible solutions}
\label{subsec:no-reducibles}

We  now turn to finding conditions which will ensure that the moduli
space contains no reducible solutions. We shall say that a class $c
\in H^{2}(X;\R)$ is \emph{integral} if it in the image of
$H^{2}(X;\Z)$ (or equivalently if it has integer pairing with every
class in $H_{2}(X;\Z)$). We shall use $\cH^{-}$ to denote the
space of real anti-self-dual harmonic $2$-forms, which we may regard
as a linear subspace of $H^{2}(X;\R)$.

\begin{proposition}
    If $M^{w}_{\kappa}$ contains a reducible solution, then there is
    an integer $n$ with $0 < n < N$ and an integral class $c \in
    H^{2}(X;\R)$ such that
    \[
           c - \frac{n}{N}c_{1}(w) \in \cH^{-}.
    \]
\end{proposition}

\begin{proof}
Suppose that $[A]$ is
a reducible solution in $M^{w}_{\kappa}$, and let $\tilde A$ be the
corresponding $U(N)$-connection in the bundle $P \to X$, inducing the
connection $\theta$ in $w$. The connection $\tilde{A}$ respects a
reduction of $P$ to $U(n_{1}) \times U(n_{2})$, with $n_{1} + n_{2} =
N$: we have
\[
        P \supset P_{1} \times_{X} P_{2},
\]
where $P_{i}$ is a bundle with structure group $U(n_{i})$, and
$n_{i} > 0$. Let $\tilde A_{1}$ and $\tilde A_{2}$ be the connections
obtained from $\tilde A$ in the bundles $P_{1}$ and $P_{2}$.  
Let $\iota : \u(1) \to \u(N)$ be the inclusion that is the derivative
of the inclusion of $U(1)$ as the center of $U(N)$. Because $\tilde A$
induces the connection $\theta$ on $w$ and $A$ is anti-self-dual, we have
\[
      F^{+}_{\tilde A} = (1/N) \iota( F^{+}_{\theta} ).
\]
Because $\tilde A = \tilde A_{1} \oplus \tilde A_{2}$, we have
\[
     F^{+}_{\tilde A_{1}} = (1/N) \iota_{1} (F^{+}_{\theta}),
\]
where $\iota_{1} : \u(1) \to \u(n_{1})$ is the inclusion again. Let
$w_{1}$ be the determinant line bundle of $P_{1}$, and $\theta_{1}$
the connection that it obtains from $\tilde A_{1}$. From the equality
above, we obtain
\[
     F^{+}_{\theta_{1}} = (n_{1}/N)  F^{+}_{\theta}.
\]
Since $F_{\theta}$ and $F_{\theta_{1}}$ represent $-2\pi i c_{1}(w)$
and $-2\pi i c_{1}(w_{1})$ respectively, we have proved the
proposition: for the class $c$ in the statement, we take
$c_{1}(w_{1})$, and for $n$ we take $n_{1}$.
\end{proof}

\begin{definition}\label{def:coprime}
    We will say that $c_{1}(w)$ is \emph{coprime to $N$} if there is
    no $n$ with $0 < n < N$ such that $(n/N)
    c_{1}(w)$ is an integral class. This is equivalent to saying that
    there is a class in $H_{2}(X;\Z)$ whose pairing with $c_{1}(w)$ is
    coprime to $N$.
\end{definition}

As corollaries to the proposition, we have:

\begin{corollary}\label{cor:codim-b+}
  Let $w$ be a fixed line bundle on $X$, and suppose $c_{1}(w)$ is
  coprime to $N$.
  Let $\cM$ denote the set of all $C^{r}$ metrics
  on $X$, and let $\Xi \subset \cM$ denote the set of metrics for
  which there exists a $\kappa$ such that $M^{w}_{\kappa}$ contains
  reducible solutions. Then $\Xi$ is contained in a countable union of
  smooth submanifolds of $\cM$, each of which has codimension
  $b^{+}_{2}$ in $\cM$.
\end{corollary}

\begin{proof}
  If the coprime condition holds, then for each integral class $c$ and
  each $n$ less than $N$,
  the class  $c - (n/N) c_{1}(w)$ is non-zero. The set of metrics
  $g$ such that
  \[
             c - (n/N) c_{1}(w) \in \cH^{-}(g)
  \]
  is therefore a smooth submanifold $\cM_{c,n} \subset \cM$, of codimension
  $b^{+}_{2}$. (See \cite{Donaldson-Kronheimer}, section 4.3, for example.)  The set
  $\Xi$ is contained in the union of these submanifolds $\cM_{c,n}$ as
  $c$ runs through the integral classes and $n$ runs through the
  integers from $1$ to $N-1$.
\end{proof}

\begin{corollary}\label{cor:reducibleSummary}
    If $c_{1}(w)$ is coprime to $N$ and $b^{+}_{2}$ is non-zero, then for a residual set of
    Riemannian metrics $g$ on $X$, the corresponding moduli spaces
    $M^{w}_{\kappa}$ contain no reducible solutions. \qed
\end{corollary}

\section{Perturbing the equations}
\label{sec:perturb}

When $N=2$, it is known that for a residual set of Riemannian metrics,
the irreducible solutions in the corresponding moduli spaces
$M^{w}_{\kappa}$ are all regular. This is the generic metrics theorem
of Freed and Uhlenbeck \cite{Freed-Uhlenbeck}. Unfortunately, no such result is
known for larger $N$, and to achieve regularity we must perturb the
anti-self-duality equations. The holonomy perturbations which we use
follow a plan which appears also in
\cite{Donaldson-Orientations}, \cite{FloerOriginal} and \cite{Taubes-Cass} among
other places. A discussion of bubbles in the presence of such
perturbations occurs in \cite{Donaldson-Book}.

\subsection{A Banach space of perturbations}

Let $X$ be a Riemannian $4$-manifold as before, let $B$ be an embedded
ball in $X$, and let
$q : S^{1} \times B \to X$ be a smooth map with the following two
properties:
\begin{enumerate}
    \item the map $q$ is a submersion;
    \item the map $q(1,\dash)$ is the identity: $q(1,x)  = x$ for all
    $b\in B$.
\end{enumerate}
Let
\[
    \omega \in \Omega^{+}(X;\C)
\]
be a smooth complex-valued self-dual 2-form whose support is contained in
$B$. For each $x$ in $B$, let 
\[
       q_{x} : S^{1} \to X
\]
be the map given by $q_{x}(z) = q(z,x)$.
Given a smooth connection $A$ in $\cA$, let $\tilde A$ be the corresponding
connection in $P \to X$ with determinant $\theta$, and let
$\Hol_{q_{x}}(\tilde A)$ denote the holonomy of $\tilde A$ around
$q_{x}$. As $x$ varies over $B$, the holonomy
defines a section of the bundle of unitary groups $U(P)$ over
$B$. We regard $U(P)$ as a subset of the vector bundle
$\gl_{P}$, with fiber $\gl(N)$. We can then multiply the holonomy by
the $2$-form $\omega$, and extend the resulting section by zero to
define a section on all of $X$:
\[
   \omega \otimes \Hol_{q}(\tilde A)
   \in  \Omega^{+}(X; \gl_{P}).
\]
We apply the orthogonal bundle projection $\pi : \gl_{P} \to \g_{P}$
to this section, and finally obtain a smooth section
\[
        V_{q,\omega}(A) \in \Omega^{+}(X; \g_{P}).
\]
Recall that $\cA$ is defined as the $L^{2}_{l}$ completion of the
space of smooth connections.

\begin{proposition}\label{prop:HolSmooth}
    For fixed $q$ and $\omega$, the map $V_{q,\omega}$
    extends to a smooth map of Banach manifolds,
    \[
          V_{q,\omega}: \cA \to
            L^{2}_{l}(X; \Lambda^{+}\otimes \g_{P}).
    \]
    Furthermore, if we fix a reference connection $A_{0}$,
    then there are constants $K_{n}$,
    depending only on $q$ and $A_{0}$, such that the $n$-th derivative
    \[
           D^{n}V_{q,\omega}|_{A} : L^{2}_{l}(X; \Lambda^{1}\otimes
           \g_{P})^{n}
                      \to L^{2}_{l}(X; \Lambda^{+}\otimes
           \g_{P})
    \]
    satisfies
    \[
          \bigl\|  D^{n}V_{q,\omega}|_{A}(a_{1},\dots, a_{n})
         \bigr \|_{L^{2}_{l,A_{0}}}
           \le K_{n} \bigl\| \omega
           \bigr \|_{C^{l}}
           \prod_{i=1}^{n}
           \| a_{i} \|_{L^{2}_{l,A_{0}}}.
    \]
\end{proposition}

\begin{proof}
    For smooth connections $A$ and $A_{0}$, consider first the pull-backs
    $q^{*}(A)$ and $q^{*}(A_{0})$ on $S^{1}\times B$. Because $q$ is a
    submersion, the pull-back map is continuous in the $L^{2}_{l}$
    topology: there is a constant $c$ depending on $q$ and $A_{0}$, such
    that
    \[
            \| q^{*}(A) - q^{*}(A_{0}) \|_{L^{2}_{l,q^{*}(A_{0})}} \le
            c \| A - A_{0} \|_{L^{2}_{l,A_{0}}}.
    \]
    In other words, the pull-back map $q^{*}$ extends continuously to
    $\cA$.

    We wish to regard $q^{*}(A)$ as providing a family of connections
    on $S^{1}$, param\-etrized by $B$. To this end, for each $x\in B$,
    let $\cH_{x}$ denote the Hilbert space of square-integrable
    $\g_{P}$-valued $1$-forms on $S^{1}\times \{x\}$:
    \[
                   \cH_{x} = L^{2}(S^{1} \times \{x\},
                   \Lambda^{1}(S^{1}) \otimes q^{*}(\g_{P})).
    \]
    These form a Hilbert bundle $\cH \to B$, and $q^{*}(A_{0})$ supplies
    $\cH$ with a connection, $A_{0}^{\cH}$. For $A$ in $\cA$, we have just observed
    that the $1$-form $q^{*}(A) - q^{*}(A_{0})$ lies in
    $L^{2}_{l}(S^{1}\times B)$; and we therefore obtain from it, by
    restriction, a section of the bundle $\cH$ of class $L^{2}_{l}$:
    \begin{equation}\label{eq:a-section}
    \begin{gathered}
            a : x \mapsto a_{x} \\
                a \in L^{2}_{l}(B;\cH).
                \end{gathered}
    \end{equation}
    (This is just the statement that an element of $L^{2}_{l}(X_{1}
    \times X_{2})$ defines an $L^{2}_{l}$ map from $X_{1}$ to
    $L^{2}(X_{2})$.) Furthermore, we have an inequality
    \[
            \| a \|_{L^{2}_{l, A_{0}^{\cH}}} \le \| q^{*}(A) -
            q^{*}(A_{0}) \|_{L^{2}_{l,q^{*}(A_{0})}}.
    \]
     We can apply these constructions also to a unitary connection
     $\tilde A$ with determinant $\theta$, to obtain $\tilde a$ in a
     similar way.

     \begin{lemma}
        Let $A_{0}$ be a smooth connection in a bundle $P \to S^{1}$,
        and write a general connection as $A = A_{0} + a$.   Identify
        the structure group $U(N)$ with the group of automorphisms of the
        fiber of $P$ over $1 \in S^{1}$. For each $a$, let $\Hol(a)
        \in U(N)$ denote the holonomy of $A_{0}+a$ around the loop
        $S^{1}$  based at $1$.  Then $\Hol$ extends to a smooth map
        \[
                \Hol : L^{2}(S^{1}; \Lambda^{1}\otimes \g_{P}) \to
                U(N).
        \]
        Furthermore, there are constants $c_{n}$, independent of
        $A_{0}$ and $A$, such that the $n$-th derivative
        $D^{n}\Hol|_{A}$ satisfies
        \[
                \bigl \| D^{n}\Hol|_{A}(a_{1}, \dots, a_{n})
                \bigr\| \le c_{n} \prod_{i=1}^{n} \| a_{i}\|_{L^{2}}.
        \]
     \end{lemma}

     \begin{proof}
        The connection $A_{0}$ does not appear in the inequality and
        serves only to define the domain of $\Hol$. As a statement
        about $A$ and the $a_{i}$, the inequality to be proved is
        gauge-invariant. It will suffice to show that, for each $A$,
        there is a connection $c_{n}(A)$ such that the inequality
        holds for all $a_{i}$: the fact that $c_{n}$ can eventually
        be taken to be independent of $A$ will follows automatically,
        because the space of connections modulo gauge is compact in
        the case of $S^{1}$.  So we are only left with the first
        assertion: that $\Hol$ extends to a smooth map on $L^{2}$. This
        is proved in \cite{Taubes-Cass, Feehan-Leness-I}. (In fact, $L^{1}$ suffices.)
     \end{proof}

     To return to the proof of Proposition~\ref{prop:HolSmooth}, we
     apply the lemma to the family of circles parametrized by $B$. For
     each $x$ in $B$, we regard the automorphism group of $P_{x}$ as
     subset of $(\gl_{P})_{x}$, and have the holonomy map
     \[
                \Hol : \cH_{x} \to (\gl_{P})_{x}.
     \]
     We regard this as giving us a smooth bundle map
     \[
            \Hol : \cH \to (\gl_{P})|_{B}.
     \]
     This map is nonlinear on  the fibers; but the uniform bounds on
     derivatives of all orders mean that it still gives a map
     \begin{equation}\label{eq:HolHilbert}
            \Hol : L^{2}_{l, A_{0}^{\cH}}(B;\cH) \to L^{2}_{l,A_{0}}(B; \gl_{P}).
     \end{equation}
     This $\Hol$ satisfies uniform bounds its derivatives also.
     (For example, if $E$ and $F$ are Hilbert spaces and $h : E \to F$
     is a uniformly Lipschitz map, then $h$ defines a uniformly
     Lipschitz map $L^{2}(B;E) \to L^{2}(B;F)$.)

     The section $V_{q,\omega}(A)$ is obtained by applying the map
     \eqref{eq:HolHilbert} to the section $a$ from
     \eqref{eq:a-section}, then multiplying by the $C^{l}$ 2-form $\omega$ and applying a linear
projection in the bundle.
\end{proof}

Now fix once and for all a countable collection of balls $B_{\alpha}$
in $X$ and maps $q_{\alpha} :
S^{1} \times B_{\alpha} \to X$, $\alpha\in \mathbb{N}$, each
satisfying the two conditions laid out above.  Let
$K_{n,\alpha}$ be constants corresponding to
$q_{\alpha}$, as in the proposition above.  Let $C_{\alpha}$ be a
sequence of positive real numbers, defined by a diagonalization, so
that 
\[
     C_{\alpha} \ge \sup \{ \, K_{n,\alpha} \mid 0\le n\le \alpha \, \}.
\]
For each $\alpha$, let $\omega_{\alpha}$ be a self-dual complex-valued
form with support in $B_{\alpha}$, and suppose that
the sum
\[
            \sum_{\alpha} C_{\alpha} \|\omega_{\alpha}\|_{C^{l}}
\]
is convergent. Then for each $n$, the series
\begin{equation}\label{eq:alpha-series}
            \sum_{\alpha} V_{q_{\alpha},\omega_{\alpha}}
\end{equation}
converges uniformly in the $C^{n}$ topology on bounded subsets of
$\cA$. That is, for each $R$, if $B_{R}(\cA)$ denotes the Sobolev ball
of radius $R$ in $\cA$ centered at $A_{0}$,
\[
     B_{R}(\cA) = \{\, A_{0} + a \mid \| a \|_{L^{2}_{l,A_{0}}} \le R
     \, \}
\]
then the series \eqref{eq:alpha-series} converges
in the uniform $C^{n}$ topology of maps $B_{R}(\cA) \to
L^{2}_{l,A_{0}}(X:\Lambda^{+}\otimes \g_{P})$. Thus the sum of the series
defines a $C^{\infty}$ map of Banach manifolds,
\[
            V : \cA \to L^{2}_{l}(X; \Lambda^{+}\otimes\g_{P}).
\]

\begin{definition}
    Fix maps $q_{\alpha}$ and constants $C_{\alpha}$ as above.
    We define $W$ to be the Banach space consisting of all sequences
    $\bomega = (\omega_{\alpha})_{\alpha\in \mathbb{N}}$ such that the sum
    \[
                 \sum_{\alpha} C_{\alpha} \| \omega_{\alpha} \|_{C^{l}}                
    \]
    converges. For each $\bomega\in W$, we write $V_{\bomega}$ for the
    sum of the series $\sum_{\alpha} V_{q_{\alpha},\omega_{\alpha}}$,
    which is  a smooth map
    \[
            V_{\bomega}  : \cA \to  L^{2}_{l}(X;\Lambda^{+}\otimes
            \g_{P}).
    \]\qed
\end{definition}

The dependence of $V_{\bomega}$ on $\bomega$ is linear, and the map
            $(\bomega, A) \mapsto V_{\bomega}(A)$
is a smooth map of Banach manifolds,
          $  W \times \cA \to L^{2}_{l}(X;\Lambda^{+}\otimes
            \g_{P})$. Given $\bomega$ in $W$, we define the perturbed
            anti-self-duality equation to be the equation
\begin{equation}\label{eq:perturbedASD}
            F^{+}_{A} + V_{\bomega}(A) = 0
\end{equation}            
for $A \in \cA$. Note that the left-hand side lies in
$L^{2}_{l-1}(X;\Lambda^{+}\otimes \g_{P})$. The equation is
gauge-invariant, and we define the perturbed moduli space to be the
quotient of the set of solutions by the gauge group:
\[
            M^{w}_{\kappa,\bomega}(X) = \{ \, [A] \in
            \B^{w}_{\kappa}(X) \mid \text{equation
            \eqref{eq:perturbedASD} holds}\,\}.
\]

\subsection{Regularity and compactness for the perturbed equations}

If $A$ is a solution of the perturbed equation
\eqref{eq:perturbedASD}, we can consider the linearization of the
equation at $A$, as a map
\begin{equation}\label{eq:derivative}
        L^{2}_{l}(X; \Lambda^{1}\otimes \g_{P}) \to  L^{2}_{l-1}(X;
        \Lambda^{+}\otimes \g_{P}) .
\end{equation}
The derivative of the perturbation $V_{\bomega}$ defines a bounded
operator from $L^{2}_{l}$ to $L^{2}_{l}$; so in the topologies of
\eqref{eq:derivative}, the derivative of $V_{\bomega}$ is a compact
operator. When we restrict the perturbed equations to the Coulomb slice through
$A$, we therefore have a smooth Fredholm map, just as in the
unperturbed case, and we have a Kuranishi theory for the perturbed
moduli space. We write the Kuranishi complex at a solution $A$ as
    \[
        L^{2}_{l+1}(X;\Lambda^{0}\otimes \g_{P})
        \stackrel{d_{A}}{\longrightarrow}
        L^{2}_{l}(X;\Lambda^{1}\otimes \g_{P})
        \stackrel{d^{+}_{A,\bomega}}{\longrightarrow}
        L^{2}_{l-1}(X;\Lambda^{+}\otimes \g_{P}) ,
\]
where $d^{+}_{A,\bomega}$ is the linearization of the perturbed
equation,
\[
       d^{+}_{A,\bomega} = d^{+}_{A} + DV_{\bomega}|_{A}.
\]
We denote the homology groups of this complex by $H^{0}_{A}$,
$H^{1}_{A,\bomega}$ and $H^{2}_{A,\bomega}$. If $A$ is irreducible and
$H^{2}_{A,\bomega}$ is zero, then the moduli space
$M^{w}_{\k,\bomega}$ is smooth near $[A]$, with tangent space
$H^{1}_{A,\bomega}$.

Understanding the Uhlenbeck compactness theorem for the perturbed
equations is  more delicate than the Fredholm theory:
if the curvature of a sequence of
connections $A_{n}$ concentrates at a point, then effect is seen in
$V_{\bomega}(A_{n})$ throughout the manifold, because the
perturbations are non-local.  To obtain a compactness result, we must
work with weaker norms than $L^{2}_{l}$.

To this end, we start by observing that
Proposition~\ref{prop:HolSmooth} continues to hold in the $L^{p}_{k}$
norms, for any $p \ge 1$ and $k\ge 0$.  That is to say, if we let
$\cA^{p}_{k}$ denote the $L^{p}_{k}$ completion of the space of smooth
connections, the $V_{q,\omega}$ defines a smooth map
\[
          V_{q,\omega} :  \cA^{p}_{k} \to L^{p}_{k}(X; \Lambda^{+}\otimes \g_{P}),
\]
and the derivatives of this map are uniformly bounded in terms of
constants which depend only on $q$ and the $C^{k}$ norm of $\omega$.
As long as the constants $C_{\alpha}$ are chosen appropriately, we can
define the Banach space of perturbations $W$ as before, and for
$\bomega\in W$ we have again a smooth map $V_{\bomega}$
\[
          V_{\bomega} :  \cA^{p}_{k} \to L^{p}_{k}(X; \Lambda^{+}\otimes \g_{P}),
\]
as long as $k\le l$.  By a diagonalization argument, we may suppose if
we wish that the constants $C_{\alpha}$ are chosen so that
$V_{\bomega}$ is smooth for \emph{all} $p$ and all $k\le l$. We can
also note that the constants $C_{\alpha}$ are inevitably bounded
below, and this implies that $V_{\bomega}(A)$ satisfies a uniform
$L^{\infty}$ bound, independent of $A$.

Consider now a sequence $[A_{n}]$ in $M^{w}_{\kappa,\bomega}$.  As just
observed, the equations imply that the self-dual curvatures
$F^{+}_{A_{n}}$ are uniformly bounded in $L^{\infty}$, and $F_{A_{n}}$
is then uniformly bounded in $L^{2}(X)$ because of the Chern-Weil
formula.  After passing to a subsequence, we may therefore assume as
usual that there is a finite set of points $\bx$ and a cover of
$X\setminus \bx$ by metric balls $\Omega_{i}$, such that
\[
        \int_{\Omega_{i}} | F_{A_{n}} |^{2} \le \epsilon
\]
for all $n$ and $i$. Here $\epsilon$ may be specified in advance to be
smaller than the constant that appears in Uhlenbeck's gauge-fixing
theorem.  We may therefore put $A_{n}$ in Coulomb gauge relative to a
trivial connection on the ball $\Omega_{i}$. For all $p\ge2$, we have a bound on
$F^{+}_{A_{n}}$ in $L^{p}$, and the connection form $A_{n}$
in Coulomb gauge on $\Omega_{i}$ is therefore bounded in $L^{p}_{1}$
by a constant depending only on $i$.

By patching gauge transformations, we arrive at the following
situation. There are gauge transformations $g_{n}$ on $X$
of class $L^{2}_{l+1,\mathrm{loc}}$, such that the sequence
$g_{n}(A_{n})$ is
bounded in $L^{p}_{1}(K)$ for any compact subset $K$ in $X\setminus
\bx$.  We rename $g_{n}(A_{n})$ as $A_{n}$.

\begin{lemma}
    Under the assumptions above, the sequence $V_{\bomega}(A_{n})$ has
    a subsequence that is Cauchy in $L^{p}(X)$.
\end{lemma}

\begin{proof}
    For $r  > 0$, let $K_{r}\subset X\setminus \bx$ be the complement
    of the open metric balls around of radius $r$ around the points of
    $\bx$. The uniform bound on the $L^{p}_{1}$ norm of $A_{n}$ on
    $K_{r}$ means that we can pass to a subsequence which converges in
    $L^{p}(K_{r})$. By  diagonalization, we can arrange the
    subsequence so that
    convergence occurs in $L^{p}(K_{r})$ for all $r > 0$.  Rename this
    subsequence as $A_{n}$. We claim that $V_{\bomega}(A_{n})$ is
    Cauchy in $L^{p}(X)$.  To see this, let $\epsilon > 0$ be given.
    For any $\beta \in \mathbb{N}$, let $V_{\bomega(\beta)}$ be the
    partial sum
    \[
               V_{\bomega(\beta)} = \sum_{\alpha=1}^{\beta}
               V_{q_{\alpha},\omega_{\alpha}}.
    \]
     Choose $\beta$ so that
     \[
                \sum_{\alpha=\beta+1}^{\infty} \| \omega_{\alpha}
                \|_{C^{0}} \le \epsilon.
     \]
     It follows that
     \[
              \| V_{\bomega}(A) -  V_{\bomega(\beta)}(A) \|_{L^{p}}
              \le \mathrm{Vol}(X)^{1/p} \epsilon
     \]
     for all $A$.  For each $r$, let $Z_{r,\alpha}\subset B_{\alpha}$
     be the set
     \[
                Z_{r,\alpha} = \{ \, x\in B_{\alpha} \mid
               \text{ $q_{\alpha}(x\times S^{1})$ is not contained in
               $K_{r}$}\, \},
     \]
     and let $Z_{r} =\cup_{\alpha\le \beta} Z_{r,\alpha}$.  The volume
     of $Z_{r}$ goes to zero as $r$ goes to zero; and because we have
     a uniform $L^{\infty}$ bound on $V_{\bomega}(A)$, we can find
     $r_{0}$ sufficiently small that
     \[
                 \| V_{\bomega(\beta)}(A) |_{Z_{r_{0}}}
                 \|_{L^{p}(Z_{r_{0}})} \le \epsilon
     \]
     for all $A$.  The $L^{p}$ convergence of $A_{n}$ on $K_{r_{0}}$
     implies that $V_{\bomega(\beta)}(A_{n})$ is Cauchy in $X\setminus
     Z_{r_{0}}$. So we can find $n_{0}$ such that, for all $n_{1}$,
     $n_{2}$ greater than $n_{0}$, we have
     \[
                 \| V_{\bomega(\beta)}(A_{n_{1}})
                 -V_{\bomega(\beta)}(A_{n_{2}})
                 \|_{L^{p}(X\setminus Z_{r_{0}})} \le \epsilon.
     \]
     Adding up the inequalities, we have
     \[
                 \| V_{\bomega}(A_{n_{1}})
                 -V_{\bomega}(A_{n_{2}}) 
                 \|_{L^{p}(X)} \le \bigl(3 + 2 \mathrm{Vol}(X)^{1/p} \bigr)\epsilon
     \]
     for all $n_{1}$, $n_{2}$ greater than $n_{0}$.  
\end{proof}

Bootstrapping and the removability of singularities now leads us to
the following version of Uhlenbeck's theorem for our situation.

\begin{proposition}\label{prop:UhlenbeckPert}
Let $A_{n}$ be a sequence of connections in $P\to X$ representing
points $[A_{n}]$ in the moduli space
$M^{w}_{\bomega,\kappa}$. Then there is a point $\mathbf{x}$ in
$X^{m}/ S_{m}$ and a connection $A'$ in a bundle
$P' \to X$ representing a point of a moduli space
$M^{w}_{\bomega,\kappa-m}$ with the following property. After passing
to a subsequence, there are isomorphisms
\[
            h_{n} : P|_{X\setminus\bx} \to P'|_{X\setminus\bx},
\]
such that $(h_{n})_{*} (A_{n})$ converges to $A'$ in $L^{p}_{1}(K)$,
for all compact $K$ contained in $X\setminus \bx$. The proof of
the lemma above shows that $(h_{n})_{*}(V_{\bomega}(A_{n}))$, extended
by zero to all of $X$, is Cauchy in $L^{p}(X)$ and has
$V_{\bomega}(A')$ as its limit. It follows that $A'$ solves the
perturbed equations. \qed
\end{proposition}

Note that, in an Uhlenbeck limit involving
bubbles, we cannot expect anything better than $L^{p}$ convergence of
the curvatures on compact subsets disjoint $\bx$, nor anything better
than $L^{p}_{1}$ convergence of the connection forms. This is in contrast
to the unperturbed case, where the convergence can be taken to be
$C^{\infty}$.

\subsection{Transversality}

Let $\cA^{*}$ again denote the irreducible connections of class
$L^{2}_{l}$. Let $W$ be the Banach space parametrizing the
perturbations $\bomega$, and consider the perturbed equations
\eqref{eq:perturbedASD} for an irreducible connection $A$ and
perturbation $\bomega$ as the vanishing of a map
\begin{equation}\label{eq:Fmap}
            \cF : W \times \cA^{*} \to
            L^{2}_{l-1}(X;\Lambda^{+}\otimes \g_{P}),            
\end{equation}
given by
\[
            \cF(\bomega, A) = F^{+}_{A} + V_{\bomega}(A).
\]
We want this map to be transverse to zero. To achieve this, we need
the family of balls $B_{\alpha}$ and maps $q_{\alpha}$ to be
sufficiently large:

\begin{condition}\label{cond:dense}
    We require the balls $B_{\alpha}$ and the submersions $q_{\alpha}: S^{1}
    \times B_{\alpha} \to X$ to satisfy the following additional
    condition: for every $x$ in $X$, the maps
    \[
                \{ \,q_{\alpha}|_{ S^{1}\times \{x\} }\mid
                \alpha\in\mathbb{N}, x \in \mathrm{int}(B_{\alpha}) \, \}
    \]
    should be $C^{1}$-dense in the space of smooth loops $q:S^{1} \to X$ based at
    $x$.
\end{condition}

\begin{lemma}\label{lem:transverse}
    Suppose that the balls $B_{\alpha}$ and maps $q_{\alpha}$ are
    chosen to satisfy the condition above. Then the map $\cF$ is
    transverse to zero.
\end{lemma}

\begin{proof}
    Let $(\bomega,A)$ be in the zero set of $\cF$.
    The derivative of $\cF$ at $(\bomega,A)$ is the map
    \[
            L : (\bnu, a) \mapsto d^{+}_{A,\bomega} a + V_{\bnu}(A).
    \]
    We shall consider the restriction of $L$ to the Coulomb slice: we
    write
    \[
            L'_{l} : W \times K_{l} \to L^{2}_{l-1}(X;\Lambda^{+}\otimes
            \g_{P}),
    \]
    where $K_{l}$ is the Coulomb slice at $A$ in the $L^{2}_{l}$
    topology. It will suffice to prove that $L'_{1}$ is surjective,
    because the regularity for the operator $d^{+}_{A}$ allows us to
    use a bootstrapping argument to show that, if $L'_{1}(\bnu, a)$
    belongs to $L^{2}_{l-1}$, then $a$ belongs to $L^{2}_{l}$.
    
    Our Fredholm theory already tells us that the cokernel of
    $L'_{1}$ has
    finite dimension. So we need only show the image of $L'_{1}$ is
    dense in $L^{2}$.  We shall show that, for a given irreducible
    $A$, the map $\bnu \mapsto V_{\bnu}(A)$ has dense image in
    $L^{2}$. 

    In order to use the hypothesis that $A$ is irreducible, we recall
    that if $\rho : H \to U(N)$ is a unitary representation of a
    group $H$, then the hypothesis of irreducibility of $H$ means that
    the complex span of the linear transformations $\rho(h)$ is all of
    $\GL(N,\C)$.  So if $x$ is any point of $X$, we an find loops
    $\gamma_{i}$ based at $x$ such that the holonomies
    $\Hol_{\gamma_{i}}(A)$, regarded as elements of $(\gl_{P})_{x}$,
   are a spanning set.

   Because of Condition~\ref{cond:dense}, we can therefore achieve the following.
   Given any $x$ in $X$, we can find an $r>0$ and $N^2$
   elements $\alpha_{i}$ ($i = 1,\dots,N^{2}$),  such that
   $B_{r_{0}}(x) \subset B_{\alpha_{i}}$ for all $i$ and such that the
   holonomies $\Hol_{q_{\alpha_{i},x}}(\tilde A)$ span $(\gl_{P})_{x}$.
   Decreasing $r$ if necessary, we may assume further that the
   sections $e_{i} = \Hol_{q_{\alpha_{i}}}(\tilde A)$ of the bundle $\gl_{P}$
   are a basis for the fiber at every point of the ball $B_{r}(x)$.
   Here we have used the continuity of these sections, which follows
   from our assumption that $L^{2}_{l}$ is contained in $C^{0}$.

   If $\omega_{i}$ are any complex-valued self-dual forms of class
   $C^{l}$  supported in $B_{r}(x)$, then the section $\sum
   \pi(\omega_{i}\otimes e_{i})$ of $\lambda^{+} \otimes \g_{P}$ can be realized
   as $V_{\bnu}(A)$, by taking $\nu_{\alpha_{i}} = \omega_{i}$. (Here
   $\pi$ is again the projection $\gl_{P} \to \g_{P}$.)  As the
   $\omega_{i}$ vary, we obtain a dense subset of the continuous
   sections supported in the ball $B_{r}(x)$.
   
    Because $x$ is arbitrary and $X$ is compact, this is enough to
    show that the sections $V_{\bnu}(A)$ are dense in $C^{0}$ as
    $\bnu$ runs through $W$. They are therefore dense in $L^{2}$ also.
\end{proof}

\begin{corollary}\label{cor:corResidual}
    Suppose that Condition~\ref{cond:dense} holds. Then
    there is a residual subset of $W$ such that for all $\bomega$ in
    this set, and all $w$ and $\kappa$, the irreducible part of the
    moduli space $M^{w}_{\kappa,\bomega}$ is regular, and therefore a
    smooth manifold of dimension $d$ given by \eqref{eq:index}.
\end{corollary}

\begin{proof}
    As usual, the zero set of $\cF$ is a $C^{\infty}$ Banach manifold
    by the implicit function theorem. The projection of $\cF^{-1}(0)$ to $W$
    is a smooth Fredholm map of index $d$ between separable Banach
    manifolds, and the fiber over $\bomega$ is the irreducible part of
    $M^{w}_{\kappa,\bomega}$; so the result follows from the
    Sard-Smale theorem.
\end{proof}

\subsection{Reducible solutions of the perturbed equations}

When $b^{+}_{2}$ is non-zero, we would like to arrange that $M^{w}_{\kappa,\bomega}$ contains no
reducible solutions. For the unperturbed case, the relevant statement
is 
Corollary~\ref{cor:reducibleSummary} above; but the discussion that
led to this corollary breaks down when the equations are perturbed.
Using the compactness theorem however, we can obtain a result when the
perturbation $\bomega$ is small enough.

\begin{proposition}\label{prop:smallPert}
    Let $\kappa_{0}$ be given, and let $g$ be a metric on $X$ with the property that the unperturbed
    moduli spaces $M^{w}_{\kappa}$ contain no reducible solutions for
    any $\kappa \le \kappa_{0}$. Then there exists $\epsilon > 0$ such
    that for all $\bomega\in W$ with $\| \bomega \|_{W} \le \epsilon$,
    the perturbed moduli spaces $M^{w}_{\kappa,\bomega}$ also contain
    no reducibles.
\end{proposition}

\begin{proof}
    Suppose the contrary, and let $[A_{n}]$ be a reducible connection
    in $M^{w}_{\kappa,\bomega_{n}}$ for some sequence $\bomega_{n}$
    in $W$ with $\|\bomega_{n}\|_{W} \to 0$ as $n\to\infty$.  The
    proof of Uhlenbeck compactness for the perturbed equations extends
    to this case, so there is a weak limit $([A], \bx)$ in
    $M^{w}_{\kappa'}$ for some $\kappa' \le \kappa$. The condition of
    being reducible is a closed condition, so $A$ is reducible and we
    have a contradiction.
\end{proof}

If we combine this proposition with the transversality results, we
obtain the following.

\begin{corollary}
    Suppose that $b^{+}_{2}$ is non-zero and $c_{1}(P)$ is coprime to
    $N$. Then
    given any $\kappa_{0}$, we can find a Riemannian metric $g$ and perturbation
    $\bomega$ such that the moduli spaces $M^{w}_{\kappa,\bomega}$ on
    $(X,g)$ are regular for all $\kappa \le \kappa_{0}$ and contain no
    reducibles. \qed
\end{corollary}

\section{Orienting the moduli spaces}

\subsection{Conventions defining orientations}
\label{subsec:OrientationConventions}

On the irreducible part of the configuration space $\B^{w}_{\kappa}$, there is a line bundle
$\Lambda$, obtained as the quotient by $\G$ of the line bundle on
$\cA$ given by the determinant of the family of operators
\begin{equation}\label{eq:D-op}
            D_{A} = d^{*}_{A } \oplus d^{+}_{A,\bomega}
\end{equation}
acting on $\g_{P}$-valued forms.
A trivialization of this line bundle provides an orientation of the
moduli space $M^{w}_{\kappa,\bomega}$ at all points where the moduli
space is regular. The space of perturbations $\bomega$ is
contractible, so in studying trivializations of $\Lambda$, we may as
well take $\bomega=0$.

The first result is that $\Lambda$ is indeed trivial: a proof for
$U(N)$ bundles is given in \cite[\S 5.4.3]{Donaldson-Kronheimer}. Our moduli spaces are
therefore orientable manifolds.  The trickier task is to isolate what
choices need to be made to in order
to specify a particular orientation, and to determine how the
orientation depends on the choices made. For this task, we follow
\cite{Donaldson-Orientations}.

First, we orient the Lie algebra
$\su(N)$. Let $\mathfrak{s}$ be the subalgebra of $\sl(N)$ consisting
of the traceless upper triangular matrices whose diagonal entries are
imaginary. The orthogonal projection $\mathfrak{s} \to \su(N)$ is a linear
isomorphism, so we can orient $\su(N)$ by specifying an orientation
of $\mathfrak{s}$.
We write $\mathfrak{s}$ as the direct sum of the diagonal
algebra $\mathfrak{h}$ and a complex vector space, the strictly upper
triangular matrices. We take the complex orientation on the latter; so
an orientation of $\mathfrak{h}$ now determines an orientation of
$\su(N)$. To orient $\mathfrak{h}$, we write $\delta_{n}$ for the
matrix with a single $1$ at the $n$-th spot on the diagonal, and we
take the ordered basis
$i(\delta_{n} - \delta_{n+1})$, $n=1,\dots,N-1$, to be an oriented
basis of $\mathfrak{h}$.

Next we choose a homology orientation for the 4-manifold $X$: that is,
we choose an orientation $o_{X}$ of the determinant line of the
operator
\[
          D =   d^{*} \oplus d^{+} : \Omega^{1}(X) \to \Omega^{0}(X)
            \oplus \Omega^{+}(X).
\]
In the case that $P \to X$ is trivial and $A$ is the trivial
connection, the operator $D_{A}$ is the obtained from $D$ by tensor
product with the vector space $\su(N)$. So
\[
       \begin{aligned}
        \ker(D_{A}) &= \ker(D)\otimes \su(N) \\
        \coker(D_{A}) &= \coker(D)\otimes \su(N).
       \end{aligned}
\]
In general, if $V$ and $W$ are oriented vector spaces,
we orient a tensor product $V\otimes W$ by choosing oriented ordered bases
$(v_{1},  \dots, v_{d})$ and $(w_{1},  \dots, w_{e})$ and taking as an
oriented basis for the tensor product the basis
\[
        (v_{1}\otimes w_{1}, \dots, v_{1}\otimes w_{e}, v_{2}\otimes
        w_{1} ,\dots).
\]
If the dimension of $W$ is even, this orientation of $V\otimes W$ is
independent of the chosen orientation of $V$. Note also that if $W$ is
complex and has the complex orientation, then this convention equips
$V\otimes W$ with its complex orientation also. Using this convention,
our homology orientation $o_{X}$ and our preferred orientation of
$\su(N)$ together determine an orientation of the determinant line of
$D_{A}$ at the trivial connection. This trivializes the line bundle
$\Lambda$ in the case of the trivial bundle $P$.
Note that if $N$ is odd, then the
dimension of $\su(N)$ is even; and in this case the orientation of
$\Lambda$ is independent of the choice of homology orientation
$o_{X}$.

Moving away from the trivial bundle, we now consider a hermitian
rank-$N$ vector bundle $E\to X$ decomposed as
\begin{equation}\label{eq:Edecomp}
        E = L_{1} \oplus \cdots \oplus L_{N},
\end{equation}
and we take $P$ to be the corresponding principal $U(N)$ bundle.  In
line with our use of the upper triangular matrices $\mathfrak{s}$, we
can decompose the
associated bundle $\g_{P}$ as the sum of a trivial bundle with fiber
$\mathfrak{h}$ and a complex vector bundle
\[
            \bigoplus_{j > i} \Hom(L_{j}, L_{i}).
\]
Using our preferred orientation of $\mathfrak{h}$ and a chosen
homology orientation $o_{X}$, we obtain in this way an orientation
\[
            O(L_{1}, \dots, L_{N}, o_{X})
\]
for the determinant line $\Lambda \to \B_{P}$. In the case that all
the $L_{n}$ are trivial, this reproduces the same orientation as we
considered in the previous paragraph. Again, this orientation is
independent of $o_{X}$ when $N$ is odd. As a special case, given a
line bundle $w \to X$, we define $O(w, o_{X})$ to be the orientation
just defined, in the case that $L_{1} = w$ and $L_{n}$ is trivial for
$n > 1$. The corresponding bundle $P$ has $\kappa= \kappa_{0}$ given by
$\kappa_{0} = - ((N-1)/(2N))c_{1}^{2}(w)$.  In this way, for each
$w$,
we orient the determinant line over $\B^{w}_{\kappa_{0}}$.

For different values of $\kappa$, the determinant lines over
$\B^{w}_{\kappa}$ and $\B^{w}_{\kappa'}$ can be compared by the
``addition of instantons'', as in \cite{Donaldson-Orientations}.  On
$S^{4}$, there is a standard $\SU(N)$ instanton with $\kappa=1$,
obtained from the standard $\SU(2)$ instanton by the standard
inclusion of $\SU(2)$ in $\SU(N)$; using a concentrated instanton and
excision, an orientation of the determinant line over
$\B^{w}_{\kappa}$ determines an orientation over $\B^{w}_{\kappa+1}$.
A precise convention can be fixed by requiring that this determination
respects the complex orientations in the case of a K\"ahler manifold
(see below).  In this way, the orientation $O(w,o_{X})$ over
$\B^{w}_{\kappa_{0}}$ determines an orientation over $\B^{w}_{\kappa}$
for all $\kappa$.

To summarize, given a homology orientation $o_{X}$ for the
$4$-manifold, we have fixed standard orientations for the smooth parts
of all the moduli spaces $M^{w}_{\kappa}$. In the case that $N$ is
odd, the choice of $o_{X}$ is immaterial.

\subsection{Comparing orientations}

Suppose $E$ is a hermitian rank-$N$ bundle and $E' = E \otimes v$ is
another, where $v$ is a complex line bundle. The top exterior powers
are related by $w' = w + Nv$. (We use additive notation here for a
tensor product of line bundles.) The configuration spaces $\B_{P}$ and
$\B_{P'}$ for the corresponding principal bundles $P$ and $P'$ are
canonically identified, as are the moduli spaces, because the adjoint
bundles $\adP$ and $\adPp$ are the same. While this provides an
identification of two moduli spaces $M^{w}_{\kappa}$ and
$M^{w+Nv}_{\kappa}$, it is not necessarily the case that this
identification respects the orientations of their smooth parts given
by $O(w, o_{X})$ and $O(w+Nv, o_{X})$ as in the previous subsection.

We follow the proof of \cite[Proposition~3.25]{Donaldson-Orientations}
closely to compare the two. The main step is to consider the case that
$X$ is K\"ahler, and compare the orientation $O(w, o_{X})$ to the
complex orientation, where $o_{X}$ is the preferred homology
orientation of the K\"ahler manifold defined in
\cite{Donaldson-Orientations}. For a bundle $E$ decomposed as
\eqref{eq:Edecomp}, the operator
\[
    D_{A} :\Omega^{1}\Bigl(\mathfrak{h} \oplus \bigoplus_{j > i}
    \Hom(L_{j}, L_{i})\Bigr)
    \longrightarrow
   \left(\Omega^{0}\oplus \Omega^{+}\right)\Bigl(\mathfrak{h} \oplus \bigoplus_{j > i}
    \Hom(L_{j}, L_{i})\Bigr)
\]
decomposes as a sum of operators $(-\bar\partial^{*} \oplus
\bar\partial)$ acting on bundles
\[
     \begin{aligned}
        &  \Omega^{0,1}
        \otimes_{\R} \mathfrak{h} \longrightarrow
        \Bigl((\Omega^{0})^{\C} \oplus \Omega^{0,2}\Bigr)
        \otimes_{\R} \mathfrak{h}, &\\
        &  \Omega^{0,1}
        \otimes_{\C} \Hom(L_{j},L_{i}) \longrightarrow
        \Bigl((\Omega^{0})^{\C} \oplus \Omega^{0,2}\Bigr)
        \otimes_{\C} \Hom(L_{j},L_{i}),&\qquad(j > i),\\
        &  \Omega^{0,1}
        \otimes_{\C} \Hom(L_{i},L_{j}) \longrightarrow
        \Bigl((\Omega^{0})^{\C} \oplus \Omega^{0,2}\Bigr)
        \otimes_{\C} \Hom(L_{i},L_{j}),&\qquad(j > i).
     \end{aligned}
\]
We treat these three summands in turn. In each case, the index of the
operator has a complex orientation (coming from the complex structure
of $\Omega^{0,i}$) which we wish to compare to our preferred
orientation.
In the first summand, one comparison that arises is between the
complex orientation of $V\otimes_{\R}\mathfrak{h}$ and the tensor
product orientation of
$V\otimes_{\C} \mathfrak{h}$ defined in the
previous subsection, when $V$ is a complex vector space. These two differ by the sign
\[
                (-1)^{(\dim_{\C}V) (N-1)(N-2)/2},
\]
irrespective of our convention for orienting $\mathfrak{h}$. This
consideration contributes a sign
\[
                (-1)^{(\ind(-\bar\partial^{*} \oplus
                              \bar\partial)) (N-1)(N-2)/2},
\]
to our comparison.  As in \cite{Donaldson-Orientations}, there is a
further consideration in the first term, in that the preferred
homology orientation of $X$ uses the opposite orientation  to
that given by the complex structure of $(-\bar\partial^{*} \oplus
\bar\partial)$, and this contributes another $(\ind(-\bar\partial^{*} \oplus
\bar\partial)) (N-1)$. Altogether, the first summand contributes a
sign
\begin{equation}\label{eq:summand1}
                (-1)^{(\ind(-\bar\partial^{*} \oplus
                              \bar\partial)) N(N-1)/2}.
\end{equation}

In the second summand, our convention for defining
$O(L_{1},\dots,L_{N}, o_{X})$ uses the complex structure of
$\Hom(L_{j},L_{i})$, and this is the same complex structure on the
index as arises from the complex structure of $\Omega^{0,i}$. So there
is no correction here. In the third summand, our complex structure is
opposite to the one inherited from $\Omega^{0,i}$. We therefore get a
correction
\begin{equation}\label{eq:summand3}
         \prod_{j>i}   (-1)^{\ind(-\bar\partial^{*} \oplus
\bar\partial)_{\Hom(L_{i},L_{j})}}.
\end{equation}

We now specialize to the case that $L_{1}=w$ and the other $L_{n}$ are
trivial, which is the case that we used to define $O(w,o_{X})$ above.
The bundles $\Hom(L_{i},L_{j})$ with $j>i$ contribute $N-1$ copies of
$\bar w$ and $(N-2)(N-1)/2$ trivial bundles. The formula
\eqref{eq:summand3}
above (the contribution of the third summand) thus simplifies to
\[
         (-1)^{\ind(-\bar\partial^{*} \oplus
\bar\partial)_{\bar w} (N-1)+ \ind(-\bar\partial^{*} \oplus
\bar\partial) (N-1)(N-2)/2 }.
\]
Now we add the contribution \eqref{eq:summand1} from the first
summand, and we obtain
\[
      (-1)^{(N-1)( \ind(-\bar\partial^{*} \oplus
\bar\partial)_{\bar w} - \ind(-\bar\partial^{*} \oplus
\bar\partial),}
\]
which is equal to
\[
       (-1)^{(N-1)( w\cdot w + K\cdot w)/2},
\]
where $K$ is the canonical class of $X$. Thus we have proved:

\begin{proposition}\label{prop:complexcompare}
    If $X$ is K\"ahler and $o_{X}$ is the preferred homology
    orientation of the K\"ahler manifold, then the distinguished
    orientation $O(w, o_{X})$ of the moduli space $M^{w}_{\kappa}$
    coincides with its complex orientation if $N$ is odd. If $N$ is
    even, then the two orientations compare by the sign
    $(-1)^{( w\cdot w + K\cdot w)/2}$. \qed
\end{proposition}

Finally, we can use Proposition~\ref{prop:complexcompare} to compare
$O(w,o_{X})$ with $O(w+Nv, o_{X})$ on a K\"ahler manifold, using the
fact that the complex orientations of $M^{w}_{\kappa}$ and
$M^{w+Nv}_{\kappa}$ coincide. We consider the case that $N$ is even.
If we set $w' = w + Nv$, then
\[
    \bigl( ( w\cdot w + K\cdot w) - ( w'\cdot w' + K\cdot w')
    \bigr)\bigm/2     =  (N/2)v\cdot v, \pmod 2.
\]
Thus we have:

\begin{proposition}\label{prop:liftcompare}
    The two orientations $O(w , o_{X})$ and $O(w + Nv, o_{X})$ of the
    moduli space $M^{w}_{\kappa} = M^{w+Nv}_{\kappa}$ are the same
    when $N$ is odd or when $N$ is zero mod $4$. When $N$ is $2$
    mod $4$, they compare with the sign
    $(-1)^{v\cdot v}$. 
\end{proposition}

\begin{proof}
    We have just dealt with the case that $X$ is K\"ahler. The general
    case is reduced to this one by excision, as in
    \cite{Donaldson-Orientations}.
\end{proof}

\subsection{Dual bundles}

The map $\delta : U(N) \to U(N)$ given by $g \mapsto \bar{g}$ allows
us to associate to our principal bundle $P$ a new principal bundle
$\bar{P}$. The corresponding vector bundles $E$ and $\bar{E}$ obtained
from the defining representation of $U(N)$ are dual. To each
connection $A$ in the adjoint bundle $\g_{P}$, we obtain a connection
$\delta(A)$ in $\g_{\bar{P}}$; and this gives us a map
\[
\delta : M^{w}_{\kappa} \to M^{-w}_{\kappa}.
\]
In the case $N=2$, this is the
same as one obtains by tensoring $E$ with $\bar{w}$; but for other
$N$ it is different.

\begin{proposition}\label{prop:dualcompare}
    Let the smooth parts of $M^{w}_{\kappa}$ and $M^{-w}_{\kappa}$ be given the
    orientations $O(w, o_{X})$ and $O(-w, o_{X})$, where $o_{X}$ is a
    chosen homology-orientation for $X$. Then the dualizing map $\delta$
    respects these orientations in the case that $N$ is odd. In the
    case that $N$ is even, the map $\delta$ preserves or reverses
    orientation according to the sign $(-1)^{w\cdot w}$. \qed
\end{proposition}

\begin{proof}
    In the case that $X$ is K\"ahler, the map $\delta$ respects the
complex orientations of the moduli spaces, for it is a
complex-analytic map between them. Using
Proposition~\ref{prop:complexcompare} and excision, we can again
compare the orientations $O(w, o_{X})$ and $O(-w, o_{X)})$. 
\end{proof}

\section{Integer invariants and product formulae}

In this section, we use the gauge theory and perturbations that we
have developed for $\PSU(N)$ to define simple integer invariants of
$4$-manifolds.  As a preliminary step for our later calculations, we
also prove a product formula for such invariants, for a $4$-manifold
$X$ decomposed along a $3$-torus.

\subsection{Simple integer invariants}

Let $X$ be a smooth, compact oriented $4$-manifold, with
$b^{+}_{2}(X) \ge 2$. Fix $N$ and let $P \to X$ be a principal $U(N)$
bundle. Let $w$ be the $U(1)$ bundle $\det(P)$, and let $\kappa$ be
the instanton number. If $N$ is even, choose also a homology
orientation $o_{X}$ for $X$.
We impose the following simplifying condition.

\begin{hypothesis}
   We shall suppose that $c_{1}(w)$ is coprime to $N$ in the sense of
   Definition~\ref{def:coprime}. We shall also suppose that the formal
   dimension $d$ of the moduli space $M^{w}_{\kappa}$,
   as given by the formula \eqref{eq:index}, is zero.
\end{hypothesis}

By Corollary~\ref{cor:reducibleSummary}, we can choose a Riemannian
metric $g$ on $X$ so that the moduli spaces $M^{w}_{\kappa'}$ contain
no reducibles solutions. More specifically, we want to choose $g$ so
that there is no integer class $c$ and no $n < N$ such that $c -
(n/N)c_{1}(w)$ is represented by an anti-self-dual form.
Having chosen such a $g$, we can find
$\epsilon(g)$ as in Proposition~\ref{prop:smallPert} such that for all
$\bomega\in W$ with $\|\bomega\|_{W} \le \epsilon(g)$, the
perturbed moduli space $M^{w}_{\kappa'}$ contains no reducibles for
any $\kappa' \le \kappa$. Finally, we use
Corollary~\ref{cor:reducibleSummary} to find a $\bomega$ with $\|
\bomega\|_{W} < \epsilon(g)$ such that the moduli spaces
$M^{w}_{\kappa'}$ are smooth for all $\kappa'\le \kappa$. We call any
$(g,\bomega)$ arrived at in this way a \emph{good pair}.

Let $(g,\bomega)$ be a good pair. The moduli space
$M^{w}_{\kappa,\bomega}$ is a smooth $0$-manifold, and is compact
because the moduli spaces $M^{w}_{\kappa-m,\bomega}$ are empty, being
of negative dimension. So $M^{w}_{\kappa,\bomega}$ is a finite set of
points. Using $O(w,o_{X})$, we can orient this moduli space: the
points of the moduli space then acquire signs; and in the usual way we
can count the points, with signs, to obtain an integer.

The integer which we obtain in this way is independent of the choice
of good pair $(g,\bomega)$, by the usual cobordism argument. In a
little more detail, suppose $(g_{0},\bomega_{0})$ and
$(g_{1},\bomega_{1})$ are good pairs. Because $b^{+}_{2}$ is greater
than $1$, Corollary~\ref{cor:codim-b+} allows us to choose a smooth
path of metrics $g_{t}$ joining $g_{0}$ to $g_{1}$, such that the
moduli spaces $M^{w}_{\kappa}(X,g_{t})$ contain no reducibles. We can
then find a continuous function $epsilon(t)$ such that
$M^{w}_{\kappa,\bomega_{t}}(X,g_{t})$ contains no reducibles whenever
$\| \bomega_{t} \|_{W_{t}}\le \epsilon(t)$.  Here $W_{t}$ is the
Banach space of perturbations on $(X,g_{t})$: we can take it that the
maps $q_{\alpha}$ and balls $B_{\alpha}$ used in its definition are
independent of $t$. Let
$\mathcal{W} \to
[0,1]$ be the fiber bundle with fiber the perturbation space $W_{t}$.
Over $[0,1]$, the maps $\cF$ of \eqref{eq:Fmap} define a total map
\[
            \cF : \mathcal{W} \times \cA^{*} \to \mathcal{L},
\]
where $\mathcal{L}$ is the fiber bundle over $[0,1]$ with fiber $L^{2}_{l-1}(X;
\Lambda^{+}_{t} \otimes \g_{P})$. We can choose a section
$\bomega_{t}$ of the $\mathcal{W}$ over $[0,1]$, transverse to $\cF$,
and always smaller than $\epsilon(t)$. The family of moduli spaces
$M^{w}_{\kappa,\bomega_{t}}(X,g_{t})$ will then sweep out an oriented
$1$-dimensional cobordism between
$M^{w}_{\kappa,\bomega_{0}}(X,g_{0})$ and
$M^{w}_{\kappa,\bomega_{1}}(X,g_{1})$.

We summarize our discussion.

\begin{definition}\label{def:integerInvariant}
    Let $X$ be a closed, oriented smooth $4$-manifold with
    $b^{+}_{2}(X) \ge 2$, equipped with a homology orientation
    $o_{X}$.  Let $w$ be a line bundle with $c_{1}(w)$ coprime to
    $N$, and suppose there exists a $U(N)$ bundle $P$, with
    $\det(P)=w$, such that the corresponding moduli space
    $M^{w}_{\kappa}$ has formal dimension $0$. Then we define an
    integer $q^{w}(X)$ as the signed count of the points in the moduli
    space $M^{w}_{\kappa,\bomega}(X,g)$, where $(g,\bomega)$ is any
    choice of good pair. This integer depends only on $X$, $w$ and
    $o_{X}$ up to diffeomorphism; and if $N$ is odd it is independent
    of $o_{X}$.

    We extend the definition of $q^{w}$ by declaring it to be zero
    if there is no $P$ for which the corresponding moduli space is
    zero-dimensional.
    \qed
\end{definition}

For given $w$, there may be no corresponding $P$ for which the moduli
space $M^{w}_{\kappa}$ has formal dimension zero. Referring to the
formula \eqref{eq:index} for $d$ and the definition of $\kappa$, we see that  a necessary and
sufficient condition is that the integer
\[
      (N^{2}-1)(b^{+}_{2}(X) - b_{1}(X) + 1)
      +2(N-1) w \cdot w
\]
should be divisible by $4N$. In the case that $N$ is even, this
implies in particular that $b^{+}_{2} - b_{1}$ should be odd, and that
$w\cdot w$ has the same parity as $(b^{+}_{2} - b_{1} + 1)/2$.
From Proposition~\ref{prop:dualcompare} we
therefore derive:

\begin{proposition}
    The integer invariant $q^{w}(X)$ satisfies $q^{w}(X)= q^{-w}(X)$
    if $N$ is odd. If $N$ is even, then $q^{w}(X)$ and $q^{-w}(X)$
    differ by the sign
    \[
            (-1)^{(b^{+}_{2}(X) - b_{1}(X) + 1)/2}.
    \]
   \qed
\end{proposition}

\subsection{Cylindrical ends}

Let $X$ be an oriented $4$-manifold with boundary a connected
$3$-manifold $Y$, and let $P \to X$ be a principal $U(N)$ bundle.
Let $\Rep^{w}(Y)$ denote the \emph{representation variety} of $Y$: the
quotient of the space of flat connections in $\g_{P}\to Y$ by the
action of the group of determinant-1 gauge transformations of
$P|_{Y}$. We make the following simplifying assumption, in order to
extend the definition of the integer invariant $q^{w}$ in a
straightforward way to such manifolds with boundary.

\begin{hypothesis}\label{hyp:Rep}
    We shall suppose that $\Rep^{w}(Y)$ consists of a single point
    $\alpha = [A_{Y}]$. We shall further suppose that $\alpha$ is
    irreducible, and non-degenerate, in the sense that
    the cohomology group $H^{1}(Y;\alpha)$ with coefficients in the
    flat bundle $(\g_{P},\alpha)$ is trivial. 
\end{hypothesis}

Once a particular representative connection $A_{Y}$ on $\g_{P}|_{Y}$
is chosen, we can extend $A_{Y}$ to all of $\g_{P}$ to obtain a
connection $A_{X}$ on $X$. But there is an essential topological
choice in this process. The determinant-1 gauge
group $\G(Y)$ is not connected: its group of components is infinite
cyclic; and the gauge transformations on $Y$ that extend to $X$ are
precisely the identity component of $\G(Y)$.

Let $\tB_{Y}$ denote the quotient of the space of connections $\cA_Y$
by the identity component $\G_{1}(Y)\subset \G(Y)$. Let
$\tilde\Rep^{w}(Y)$ denote the subset of $\tB_{Y}$ consisting of flat
connections: this is an infinite set acted on transitively by the
infinite cyclic group $\G(Y)/ \G_{1}(Y)$. Let $\tilde\alpha$ be a
choice of a point in $\tilde\Rep^{w}_{Y}$. Let $A_{Y}$ be a connection
representing $\tilde\alpha$, and let $A_{X}$ be any extension of
$A_{Y}$ to all of $X$.

\medskip
Up to this point, we have made no essential use of a metric on $X$. We
now choose a Riemannian metric on $X$ that is cylindrical in a collar
of the boundary. We also suppose that $A_{X}$ is flat in the collar.
Form a cylindrical-end manifold $X^{+}$ by attaching a cylinder
$[0,\infty) \times Y$ to $X$, and extend $A_{X}$ as a flat connection
on the cylindrical part. We now introduce a space of connections on
$X^{+}$,
\[
            \cA(X^{+}; \tilde\alpha) = \{ \, A \mid A - A_{X} \in L^{2}_{l,
            A_{X}}(X^{+}; \Lambda^{1}\otimes \g_{P} \, \},
\]
a gauge group
\[
            \G(X^{ +};\tilde\alpha) = \{\, g \mid g-1 \in L^{2}_{l+1,A_{X}}(X^{+};
            \SU(P))\, \},
\]
and the quotient space $\B(X^{+};\tilde\alpha)$. The gauge group acts
freely, and the quotient space is a Banach manifold. Inside this
Banach manifold is the moduli space of anti-self-dual connections:
\[
            M^{w}(X^{+};\tilde\alpha) = \{ \, [A] \in
            \B(X^{+};\tilde\alpha)
            \mid F^{+}_{A} = 0 \, \}.
\]

We have the usual results concerning moduli spaces on cylindrical-end
manifolds. Any finite-energy anti-self-dual connection $A$ in the bundle $P \to
X^{+}$ is gauge-equivalent to a connection in
$M^{w}(X^{+};\tilde\alpha)$, for some unique $\tilde\alpha \in
\tilde\Rep^{w}(Y)$. Given $[A]$ in this moduli space, there is an
elliptic complex
\[
        L^{2}_{l+1,A_{X}}(X^{+}; \Lambda^{0}\otimes \g_{P})
        \stackrel{d_{A}}{\to}
        L^{2}_{l,A_{X}}(X^{+}; \Lambda^{1}\otimes \g_{P})
        \stackrel{d_{A}^{+}}{\to}
        L^{2}_{l-1,A_{X}}(X^{+}; \Lambda^{+}\otimes \g_{P}).
\]
We write $H^{i}_{A}$ for its cohomology groups. We always have
$H^{0}_{A}=0$, and again describe $[A]$ as regular if $H^{2}_{A}$ is
zero, in which case the moduli space is smooth near $[A]$. The
dimension of the moduli space at regular points is the index of the
operator 
\begin{equation}\label{eq:Fredholm-cylinder}
         d^{*}_{A_{X}} \oplus d^{+}_{A_{X}} :
          L^{2}_{l,A_{X}}(X^{+}; \Lambda^{1}\otimes \g_{P})
          \longrightarrow
           L^{2}_{l-1,A_{X}}(X^{+}; (\Lambda^{0}\oplus \Lambda^{+}) \otimes \g_{P}).       
\end{equation}
We write $d(X^{+};\tilde\alpha)$ for this index. We write the
transitive action of $\Z$ on $\tilde\Rep^{w}(Y)$ using the notation
$\tilde\alpha \mapsto \tilde\alpha + k$. The sign convention is
fixed so that
\[
            d(X^{+};  \tilde\alpha + k) = d(X^{+}; 
            \tilde\alpha )  + 4N k.
\]
We also have
\[
            \kappa(X^{+}; \tilde\alpha + k) = \kappa(X^{+}; \tilde\alpha) + k,
\]
where $\kappa$ is defined by the same Chern-Weil integral that would
define the characteristic class \eqref{eq:kappa-def} in the closed case:
\begin{equation}\label{eq:ChernWeil}
\begin{aligned}
 \kappa(X^{+}; \tilde\alpha) &= -\frac{1}{2N} \int_{X^{+}} p_{1}(\g_{P};A) \\
            &= \frac{1}{16N\pi^{2}} \int_{X^{+}} \mathrm{tr}\bigl(
            \mathrm{ad}(F_{A}) \wedge \mathrm{ad}(F_{A}) \bigr).
            \end{aligned}
\end{equation}
(The notation $\mathrm{ad}(F_{A})$ denotes a $2$-form with values in
$\mathrm{End}(\g_{P})$, and $\mathrm{tr}$ is the trace on
$\mathrm{End}(\g_{P})$.)

We can construct a Banach space $W$ parametrizing perturbations of the
anti-self-duality equations on $X^{+}$. In order to achieve
transversality, we again choose a family of balls and maps
$q_{\alpha}$ satisfying the density condition (Condition
\ref{cond:dense}). The constants $C_{\alpha}$ in the definition of $W$
can be chosen so that for all $\bomega\in W$, the perturbing term
$V_{\bomega}(\alpha)$ has rapid decay on the end of $X^{+}$: for some
constants $K_{j}$, and $t$ a function equal to the first coordinate on
the cylindrical end, we can require
\[
     |  \nabla^{j}_{A_{X}} V_{\bomega}(A) | \le K_{j} e^{-t}  \|\bomega
     \|_{W} ,
\]
for all $A$ in $\cA(X^{+};\tilde\alpha)$, all $\bomega$ in $W$, and
all $j \le l$. With such a condition, the term $V_{\bomega}(A)$ will
contribute a compact perturbation to the linearized equations, and we
have moduli spaces $M^{w}_{\bomega}(X^{+};\tilde\alpha)$. For a
residual set of perturbations $\bomega$ in $W$, the perturbed moduli
space $M^{w}_{\bomega}(X^{+};\tilde\alpha)$ will be regular, and
therefore a smooth manifold of dimension $d(X^{+};\tilde\alpha)$.
A pair $(g,\bomega)$ consisting of a cylindrical-end metric and
perturbation $\bomega$ will be called a \emph{good pair} if all the
moduli spaces $M^{w}_{\bomega}(X^{+};\tilde\alpha)$ are regular.

The
moduli space is orientable, and can be oriented by choosing a
trivialization of the determinant line bundle on the connected Banach
manifold $\B^{w}(X^{+};\tilde\alpha)$. We make no effort here to
define a canonical orientation for the moduli space.

The proof of the compactness theorem also adapts to the cylindrical
end case. If $[A_{n}]$ is a sequence in
$M^{w}_{\bomega}(X^{+};\tilde\alpha)$, then after passing to a
subsequence there is an Uhlenbeck limit $( [A]; \bx )$, where $\bx \in
\mathrm{Sym}^{m}(X^{+})$ and $[A] \in
M^{w}_{\bomega}(X^{+};\tilde\alpha-m')$. The main difference from the
case of a closed manifold is that we only have $m' \ge m$, rather than
equality, because some energy may be lost on the cylindrical end.

We can now define an integer invariant (with an ambiguous sign)
\begin{equation}\label{eq:withBoundary-qw}
        q^{w}(X;\alpha ) \in \Z / \{\pm 1\} 
\end{equation}
whenever Hypothesis~\ref{hyp:Rep} holds. We define it by picking a
good pair $(g,\bomega)$, choosing an
orientation for the determinant line, and then counting with signs the
points in a zero-dimensional moduli space
$M^{w}_{\bomega}(X^{+};\tilde\alpha)$ if one exists. If there is no
$\tilde\alpha$ for which the moduli space is zero-dimensional, we
define $q^{w}(X;\alpha)$ to be zero.

Note that there is no hypothesis on $b^{+}_{2}(X)$ in this
construction, because there are no reducible solutions in the moduli
spaces.

\subsection{Perturbations with compact support}

The following lemma is useful for the gluing formula in the subsection
below.

\begin{lemma}\label{lem:compactSupport}
    Let $X^{+}$ be a manifold with cylindrical end, as above, and
    suppose that $w$ satisfies Hypothesis~\ref{hyp:Rep}. Suppose
    $M^{w}_{\kappa}(X^{+};\tilde\alpha)$ has formal dimension $0$. Then we can
    find a cylindrical-end metric $g$ and a perturbation $\bomega$
    making the moduli spaces $M^{w}_{\kappa,\bomega}(X^{+};\tilde\alpha-m)$
    regular for all $m\ge 0$, with the additional condition that the
    sum \eqref{eq:alpha-series} defining the perturbation has only finitely
    many non-zero terms: that is, only finitely many of the
    $\omega_{\alpha}$ are non-zero.
\end{lemma}

\begin{proof}
    We already know that we can find $\bomega$ in $W$ such that the
    perturbed moduli spaces
    $M^{w}_{\kappa,\bomega}(X^{+};\tilde\alpha-m)$ are regular. 
    (For $m$ strictly positive, regular means that this
    moduli space is empty.) We will show that the regularity of
    these moduli spaces for all $m \ge 0$ is an open condition. This
    will suffice, because the perturbations given by finite sums are
    dense.

    Suppose then $\bomega_{n}$ is a sequence converging to $\bomega$ in
    $W$. Because of the Chern-Weil formula, we know that there exists
    $m_{0}$ such that $M^{w}_{\kappa,\bomega_{n}}(X^{+};\tilde\alpha-m)$
    is empty (and therefore regular), for all $m\ge m_{0}$ and all
    $n$.

    Suppose that $m_{0}$ is greater than $1$, and consider the moduli
    spaces
    \[M^{w}_{\kappa,\bomega_{n}}(X^{+};\tilde\alpha-m_{0}+1).\] If these
    are non-empty for infinitely many $n$, then we have a
    contradiction, arising from Uhlenbeck's compactness theorem and
    the emptiness of $M^{w}_{\kappa,\bomega}(X^{+};\tilde\alpha-m)$
    for $m$ positive. In this way, we prove inductively that
    $M^{w}_{\kappa,\bomega_{n}}(X^{+};\tilde\alpha-m)$ is empty for
    all $n\ge n_{0}$ and all $m \ge 1$.

    Now consider the moduli spaces
    $M^{w}_{\kappa,\bomega_{n}}(X^{+};\tilde\alpha)$, which have
    formal dimension $0$. Suppose these are non-regular for infinitely
    many $n$. After passing to a subsequence, we may suppose there is
    an irregular solution $[A_{n}]$ in
    $M^{w}_{\kappa,\bomega_{n}}(X^{+};\tilde\alpha)$ for all $n$, and
    that these converge in Uhlenbeck's sense. The fact that the lower
    moduli spaces are empty means that this is in fact strong
    convergence of the connections $A_{n}$ after gauge transformation.
    Irregularity is a closed condition under strong limits, so the
    limit is an irregular point $[A]$ in
    $M^{w}_{\kappa,\bomega}(X^{+};\tilde\alpha)$, which contradicts
    our hypothesis.
\end{proof}

\begin{remark}
    One should expect to prove that irregularity is a closed condition
    under Uhlenbeck limits, and so extend the lemma to
    higher-dimensional moduli spaces. Because of the non-local nature
    of our perturbations, this would require some more work.
\end{remark}

\subsection{Gluing}
\label{subsec:gluing}

Suppose $X$ is a closed $4$-manifold with $b^{+}_{2}\ge 2$, so that the integer
invariants $q^{w}(X)$ are defined. Suppose $X$ contains a connected
$3$-manifold $Y$ that separates $X$ into two manifolds with common boundary,
$X_{1}$ and $X_{2}$, and suppose that the representation variety
$\Rep^{w}(Y)$ satisfies Hypothesis~\ref{hyp:Rep}. Let $w_{i}$ be the
restriction of $w$ to $X_{i}$, so that we have integer invariants
(with ambiguous sign) $q^{w_{i}}(X_{i};\alpha)$ for $i=1,2$.

The sign-ambiguity can be partly resolved as follows. Let $\Lambda$ be
the orientation line bundle on $\B^{w}(X)$, and let $\Lambda_{i}$ be
the orientation bundles on $\B^{w_{i}}(X^{+};\tilde\alpha)$. Then
there is a preferred isomorphism
\[
            \Lambda = \Lambda_{1}\otimes \Lambda_{2}.
\]
Suppose we choose a homology orientation $o_{X}$ for $X$, and so
trivialize $\Lambda$ using $O(w, o_{X})$. Then an orientation for
$\Lambda_{1}$ determines an orientation for $\Lambda_{2}$, using the
above product rule. Thus a choice of sign for
$q^{w_{1}}(X_{1};\alpha)$ determines a choice of sign for
$q^{w_{2}}(X_{2};\alpha)$. In the next proposition, we assume that the
signs are resolved in this way.

\begin{proposition}\label{prop:Product}
    In the above situation, we have a product law,
    \[
            q^{w}(X) = q^{w_{1}}(X_{1};\alpha)
            q^{w_{2}}(X_{2};\alpha).
    \]
\end{proposition}

\begin{proof}
    First, choose perturbation $\bomega_{1}$ and $\bomega_{2}$ for
    $X_{1}^{+}$ and $X_{2}^{+}$ satisfying the condition in the
    conclusion of Lemma~\ref{lem:compactSupport}. This condition means
    that there is a compact subset $K_{i}\subset X_{i}^{+}$ such that,
    for all $A$ on $X^{+}_{i}$, the perturbing term
    $V_{\bomega_{i}}(A)$ is supported in $K_{i}$ and depends only on
    the restriction of $A$ to $K_{i}$.  Equip $X$ with a metric
    $g_{R}$ containing a long cylindrical neck $[-R,R]\times Y$, in
    the usual way. For $R$ sufficiently large, $X$ contains isometric
    copies of $K_{1}$ and $K_{2}$, so we can regard $\bomega_{1}$ and
    $\bomega_{2}$ as defining a perturbation $V_{\bomega}$ of the
    equations on $(X, g_{R})$. Since this perturbation is supported
    away from the neck region, it does not interfere with the standard
    approaches to gluing anti-self-dual connections. The conclusion is
    that, for $R$ sufficiently large, the moduli space
    $M^{w}_{\bomega}(X, g_{R})$ is regular and is the product of the two moduli
    spaces $M^{w_{i}}_{\bomega_{i}}(X^{+}_{i};\tilde\alpha)$, provided
    there exists a lift $\tilde\alpha$ such that these are
    zero-dimensional.
\end{proof}

\subsection{Other cases}
\label{subsec:othercases}

We now consider a slightly more general setting, in which the
representation variety $\Rep^{w}(Y)$ is no longer required to be a
single point. We still ask that $c_{1}(w)$ is coprime to $N$ on $Y$,
but we suppose now that $\Rep^{w}(Y)$ consists of $n$ irreducible
elements $\alpha_{i}$ ($i=1, \dots , n$), all of which are
non-degenerate.  Rather than develop a complete Floer homology theory,
we continue to make some simplifying assumptions, which we now lay
out.

If $\gamma: [0,1]\to \B^{w}(Y)$ is a path joining $\alpha_{i}$ to
$\alpha_{j}$, we can associate to $\gamma$ two quantities, both of
which are additive along composite paths. The first is the spectral
flow of the family of operators
\[
            \begin{bmatrix}
                    0 & - d^{*}_{A(t)} \\
                    -d_{A(t)} & *d_{A(t)}
            \end{bmatrix}
\]
acting on $(\Omega^{0} \oplus \Omega^{1})(Y;\g_{P})$, where $A(t)$ is
a lift of the path $\gamma$ to $\cA(Y)$. We call this integer
$d(\gamma)$. The second quantity is the Chern-Weil integral of the
same kind as \eqref{eq:ChernWeil}. We regard $A(t)$ as defining a
connection $A$ on $[0,1]\times Y$ in temporal gauge, and define
\[
            \kappa(\gamma) = 
            \frac{1}{16N\pi^{2}} \int_{[0,1]\times Y} \mathrm{tr}\bigl(
            \mathrm{ad}(F_{A}) \wedge \mathrm{ad}(F_{A}) \bigr).
\]
Modulo $\Z$, we can interpret $\kappa(\gamma)$ as the drop in the
suitably-normalized Chern-Simons invariant from $\alpha_{i}$ to
$\alpha_{j}$. The normalization is such that, for a suitably-oriented
closed loop in $\B^{w}(Y)$ representing the generator of first
homology, the value of $\kappa(\gamma)$ is $1$. For the same loop,
$d(\gamma)$ is $4N$.

If $A(t)$ is a path of connections lifting $\gamma$, and $A$ is the
corresponding $4$-dimensional connection as above, we may extend $A$
to a connection $A_{Z}$ on the infinite cylinder $Z = \R \times Y$,
just as we did for the cylindrical-end case. The connection $A_{Z}$ is
flat on both ends of th cylinder. We can then define a
space of connections
\[
            \cA = \{ \, A_{Z} + a \mid a \in L^{2}_{l,a_{Z}}(Z;
            \g_{P}) \, \}
\]
and a quotient space $\B = \B^{w}(Z;\gamma)$. Inside the configuration
space is the moduli space $M^{w}(Z;\gamma)$, and $d(\gamma)$ is its
formal dimension.

\begin{hypothesis}\label{hyp:monotone}
   Let $\gamma$ denote a path from $\alpha_{i}$ to $\alpha_{j}$ as
   above. We impose two conditions:
   \begin{enumerate}
    \item $d(\gamma)$ is not $1$ for any such $\gamma$;
    \item if $d(\gamma)\le 0$, then $\kappa(\gamma)\le 0$ also.
   \end{enumerate}
\end{hypothesis}

Assume these two conditions hold, and consider  a moduli space
$M^{w}_{\bomega}(X^{+};\tilde\alpha_{k})$ on a cylindrical-end
manifold $X^{+}$. Suppose $\bomega$ is chosen so that all moduli
spaces are regular, and that this particular moduli space is
zero-dimensional. The two conditions imply that if
$M^{w}_{\bomega}(X^{+};\tilde\alpha_{j})$ is another moduli space,
with $\kappa(X^{+};\tilde\alpha_{j}) \le \kappa(X^{+};\tilde\alpha_{k})$, then the
formal dimension of this moduli space is $-2$ or less. This moduli
space is therefore empty, and remains empty for a generic path of
metrics an perturbations. It follows that there is a well-defined
integer invariant
\[
        q^{w}(X;\alpha_{k}) \in \Z/ \pm 1
\]
counting points in this zero-dimensional moduli space. Once again, we
extend the definition by declaring it to be zero if there is no lift
$\tilde \alpha_{k}$ for which the moduli space is zero-dimensional.

The product law from Proposition~\ref{prop:Product} then extends to
this more general situation. In the situation considered in the
proposition, we now have:
\begin{equation}\label{eq:generalProduct}
            q^{w}(X) = \sum_{k=1}^{n}q^{w_{1}}(X_{1};\alpha_{k})
            q^{w_{2}}(X_{2};\alpha_{k}).
\end{equation}

\section{Calculations}

\subsection{The $K3$ surface}
\label{sec:K3}

Let $X$ be a $\mathit{K3}$ surface (with the complex orientation) and let $P
\to X$ be a $U(N)$ bundle with
\[
\begin{aligned}
               \bigl\langle c_{2}(P), [X] \bigr\rangle &= N(N^{2}-1) \\
                \bigl\langle c_{1}(P)^{2}, [X] \bigr\rangle &=
                2(N+1)^{2}(N-1).
\end{aligned}
\]
More specifically, we suppose these conditions are achieved by
having
\[
\begin{aligned}
               c_{1}(P) &=  \binom{N+1}{1}  h, \\
               c_{2}(P) &=  \binom{N+1}{2} h^{2}
\end{aligned}
\]
where $h$ is a primitive class of square $2(N-1)$. These conditions
ensure in particular that $c_{1}(P)$ is coprime to $N$.
The bundle $P$ has $\kappa = N- (1/N)$, and the formal dimension of the
corresponding moduli
space $M^{w}_{\kappa,\bomega}$ is zero.  
We choose the homology
orientation $o_{X}$ arising from the K\"ahler structure, and
Proposition~\ref{prop:complexcompare} tells us that the moduli space
has the same orientation as the one it obtains as a
complex-analytic space. In this situation, there is an integer
invariant $q^{w}(X)$ as in Definition~\ref{def:integerInvariant}, and our aim is to
calculate it.

\begin{proposition}
    The above invariant $q^{w}(X)$ for the $K3$ surface $X$ is $1$.
\end{proposition}

\begin{proof}
    We use a familiar circle of ideas to evaluate the invariant: see
    \cite{Donaldson-polynomials} for the calculation in the case
    $N=2$, on which our argument is based. We assume $N\ge 3$ for the
    rest of this proof.

    We take $X$ to be an algebraic surface embedded in $\CP^{N}$ by
    the complete linear system of a line bundle with first Chern class
    $h$. This can be done in such a way that $h$ generates the Picard
    group of $X$: all holomorphic line bundles on $X$ are powers of
    the hyperplane bundle $H$.

    For the inherited metric $g$, the corresponding moduli space
    $M^{w}_{\kappa}$ can be identified with the moduli space of poly-stable
    holomorphic bundles with the topology of $P$. The condition on
    $c_{1}(P)$ ensures that there are no reducible bundles, so
    poly-stable means stable in this case. As usual
    on a $K3$, irreducible stable bundles are regular, so the moduli
    space $M^{w}_{\kappa}$ is smooth and $0$-dimensional; and the
    moduli spaces $M^{w}_{\kappa-m}$ are empty for $m$ positive. We
    can therefore use the stable bundles to calculate the integer
    invariant. By the remark about orientations above, the integer
    $q^{w}$ is simply the number of stable bundles on $X$ with the
    correct topology.

    Mukai's argument \cite{Mukai-1984} shows that the moduli space is connected, so to
    prove the proposition we just have to exhibit one stable bundle.
    The Chern classes of $P$ have been specified so that the
    restriction to $X$ of the holomorphic tangent bundle $T\CP^{N}$ as
    the correct topology. We will show that the holomorphic vector
    bundle
    \[
            \mathcal{E} = T\CP^{N}|_{X}
    \]
    is stable, under the hypothesis that $H$ generates the Picard
    group. To do this, we recall that $\mathcal{E}$ has a description
    as a quotient
    \[
                0 \to \mathcal{O} \stackrel{s}{\to} V\otimes H
                \stackrel{q}{\to} \mathcal{E} \to 0,
    \]
    where $V$ is a vector space of dimension $N+1$. The bundle
    $\mathcal{E}$ has rank $N$ and determinant $(N+1)H$; so if it is
    unstable then there is subsheaf $\mathcal{F}$ with rank $n$ and
    determinant at least $(n+1)H$, where $n$ lies in the range $1\le n
    \le N$. We can suppose that $\mathcal{E}/\mathcal{F}$ is
    torsion-free. Let $\tilde{\mathcal{F}}$ be the inverse image of
    $\mathcal{F}$ under the quotient map $q$. Then
    $\tilde{\mathcal{F}}$ has rank $n+1$ and determinant at least
    $(n+1)H$; and the quotient of $V\otimes H$ by
    $\tilde{\mathcal{F}}$ is torsion-free.  This forces
    \[
            \tilde{\mathcal{F}} = U \otimes H,
    \]
    where $U$ is a non-trivial proper subspace of $V$. However, the
    image of the section $s$ in the above sequence is not contained in
    $U\otimes H$, for any proper subspace $U$. (Geometrically, this is
    the statement that  $X$ is not contained in any hyperplane in
    $\CP^{N}$.) Since $\tilde{\mathcal{F}}$ contains $s(\mathcal{O})$
    by construction, we have  a contradiction.
\end{proof}

\subsection{Flat connections on the 3-torus}

We consider the $3$-torus $T^{3}= \R^{3}/\Z^{3}$, and take $x$, $y$
and $z$ to be the simple closed curves which are the images of the
three coordinate axes in $\R^{3}$. We regard these as oriented curves,
defining elements of the fundamental group $\pi_{1}(T^{3}, z_{0})$,
where $z_{0}$ is the base point. We take $w\to T^{3}$ to be a line
bundle with $c_{1}(w)$ Poincar\'e dual to $z$, and we take $P$ to be a
principal $U(N)$ bundle having $w$ as determinant.

\begin{lemma}\label{lem:RepT3}
    For the $3$-manifold $T^{3} = \partial Z_{K}$, the representation
    variety $\Rep^{w}(T^{3})$ consists of $N$ points
    $\alpha_{0},\dots,\alpha_{N-1}$. Each of 
    these points is irreducible and non-degenerate.
\end{lemma}

\begin{proof}
    Let $z' \subset T^{3}$ be a parallel copy of $z$, disjoint from
    the curves $x$ and $y$ and meeting the 2-torus $T^{2}\times \{0\}$
    spanned by $x$ and $y$
    transversely in a single point. Let $\eta$ be a $2$-form
    supported in a tubular neighborhood of $z'$, representing the dual
    class to $z$. The support of $\eta$ should be disjoint from $x$,
    $y$ and $z$.  Equip the line bundle $w$ with a connection $\theta$
    whose curvature is $-(2\pi i) \eta$ and which is trivial on the
    complement of the tubular neighborhood. Use $\theta$ to trivialize
    $w$ away from $z'$. This trivialization of the determinant reduces
    the structure group of $P$ to $\SU(N)$ on the complement of the
    tubular neighborhood of $z'$.

    Using $\theta$, we identify $\Rep^{w}(T^{3})$ with the moduli space of
    flat $\SU(N)$ connections $\tilde A$ on $T^{3}\setminus z'$ such
    that the corresponding homomorphism
    \[
                \rho : \pi_{1}\bigl(T^{3}\setminus z', z_{0}\bigr) \to \SU(N)
    \]
    satisfies
    \[
                [\rho(x), \rho(y)] = e^{2\pi i /N} \mathbf{1}_{N}.
    \]
    The last condition implies that, up to similarity $\rho(x)$ and
    $\rho(y)$ can be taken to be the pair
    \[
           \begin{aligned}
            \rho(x) &= \begin{bmatrix}
                1 & 0 & 0 & \cdots& 0\\
                0 & \zeta  & 0 & \cdots&0\\
                0 & 0 & \zeta^{2} & \cdots &0\\
                  &   &        &   \ddots&0 \\
               0   &   0&   0 &  0 & \zeta^{N-1}
            \end{bmatrix}
            &
            \rho(y) &= \begin{bmatrix}
                0 & 0 & 0 & \cdots&  \pm 1\\
                1 & 0  & 0 & \cdots&0\\
                0 & 1 & 0 & \cdots &0\\
                  &   &        &   \ddots&0 \\
               0   &   0&   0 &  1 & 0
            \end{bmatrix},
           \end{aligned}
    \]
    where $\zeta= e^{2\pi i /N}$ and the sign in the top row of
    $\rho(y)$ is negative if $N$ is
    even (so that $\rho(y)$ has determinant $1$). As an element of
    $\pi_{1}(T^{3}\setminus z' , z_{0})$, the class of $z$ commutes
    with $x$ and $y$, and this tells us that
    \[
                \rho(z) = \zeta^{k} \mathbf{1}_{N}
    \]
    for some $k$ with $0\le k \le N-1$. The elements $x$, $y$ and $z$
    generate the fundamental group, so there are exactly $N$ elements
    in the representation variety, as claimed. Let $\rho_{k}$ be the
    representation with $\rho(z) = \zeta^{k} \mathbf{1}_{N}$, and let
    $\alpha_{k}$ be the corresponding point of $\Rep^{w}(T^{3})$.

    We have already remarked that the only elements commuting with
    $\rho_{k}(x)$ and $\rho_{k}(y)$ are the central elements; so each
    $\rho_{k}$ is irreducible. Each $\rho_{k}$ determines a
    representation $\bar\rho_{k}$ of $\pi_{1}(T^{3}, z_{0})$ in the
    adjoint group $\PSU(N)$; and all these are equal. To see that they
    are non-degenerate, we look at the cohomology $H^{1}(T^{3}; (\g_{P},
    \bar{\rho_{k}})$ with coefficients in the adjoint bundle. Any
    element of this group is represented by a covariant-constant
    $\g_{P}$-valued form, because the torus is flat; and there are no
    non-zero covariant-constant forms because the bundle is
    irreducible.
\end{proof}

\begin{lemma}
      For each $k$, there is a path $\gamma$ joining $\alpha_{k}$ to
      $\alpha_{k+1}$ in $\B^{w}(T^{3})$ such that  he spectral flow
      $d(\gamma)$ is $4$. For the same path, the energy
      $\kappa(\gamma)$ is $1/N$.
\end{lemma}

\begin{proof}
    The Chern-Simons invariant of $\alpha_{k}$ in $\R/\Z$ is
    calculated in \cite{Borel-Friedman-Morgan}, where it is shown to be
    $-k/N$ mod $\Z$. It follows that there is a path $\gamma_{0}$ joining
    $\alpha_{0}$ to $\alpha_{1}$ with $\kappa(\gamma_{0})=1/N$. The
    connections represented by $\alpha_{k}$ and $\alpha_{k+1}$ are
    gauge-equivalent under the larger gauge group consisting of all
    $\PSU(N)$ gauge transformations of $\adP$. It follows that we can
    lift $\gamma_{0}$ to a path of connections $A: [0,1] \to
    \cA^{w}(T^{3})$ with $A(1) =
    g(A(0))$ for some automorphism $g$ of $\adP$. Applying $g$ to this
    whole path, we obtain a path of connections $A:[1,2] \to
    \cA^{w}(T^{3})$ joining $A(2)$ to a representative $A(3)$ for some
    other point of $\Rep^{w}(T^{3})$, which must be $\alpha_{3}$ on
    the grounds of its Chern-Simons invariant.

    In this way, we construct paths $\gamma_{k}$ from $\alpha_{k}$ to
    $\alpha_{k+1}$, each of which has the same energy and the same
    spectral flow. The composite path is a loop $\delta$ based at
    $\alpha_{0}$ with $\kappa(\delta) = 1$. It follows that
    $d(\delta)$ is $4N$. So each $\gamma_{k}$ has $d(\gamma_{k}) = 4$.
\end{proof}

\subsection{Knot complements}
\label{subsec:KnotComplements}

Let $K$ be a knot in $S^{3}$, and let $M$ be the knot complement,
which we take to be a $3$-manifold with $2$-torus boundary carrying a
distinguished pair of oriented closed curves $m$ and $l$, the meridian and
longitude, on its boundary.

We write $Z_{K}$ for  the $4$-manifold $S^{1}\times M$. The boundary of
$Z_{K}$ is a $3$-torus. We identify $m$ and $l$ with curves $1\times
m$ and $1\times l$ on $\partial Z_{K}$, and write $s$ for the curve
$S^{1}\times z_{0}$ on $\partial Z_{K}$, where $z_{0}$ is the point
of intersection of $m$ and $l$, which we take as base-point in
$Z_{K}$. The three classes
\[
        [s], [m], [l] \in H_{1}(\partial Z_{K})
\]
generate the first homology of the boundary. We orient $K$, so as to
orient the longitude $l$. We suppose $m$ is oriented so that the
oriented tangent vectors of $s$, $m$ and $l$ (in that order) are an
oriented basis for the tangent space to $\partial Z_{K}$ at $z_{0}$.

The first homology of $Z_{K}$
itself is generated by $[s]$ and $[m]$, as the class $l$ is the
boundary of an oriented Seifert surface for $K$.
We take $\Sigma$ to be such a Seifert surface, regarded as a
$2$-dimensional submanifold with boundary in $(Z_{K},\partial Z_{K})$,
and we
take $w\to Z_{K}$ to be the line bundle with $c_{1}(w)$ dual to
$[\Sigma,\partial\Sigma] \in H_{2}(Z_{K},\partial Z_{K})$.

We identify $\partial Z_{K}$ with the standard $3$-torus $T^{3}$ in
such a way that $s$, $m$ and $l$ are identified with $x$, $y$ and $z$
(the curves considered in the previous subsection). From the
discussion above, we have connections $\alpha_{k}$ in
$\Rep^{w}(\partial Z_{K})$ for $k=0,\dots, N-1$, with spectral flow
$4$ along paths from each one to the next.
Because of the lemmas above, the conditions of
Hypothesis~\ref{hyp:monotone} hold for the $3$-torus $\partial
Z_{K}$, so we can consider the integer-valued
invariants
$q^{w}(Z_{K};\alpha_{k})$, as defined in section
\ref{subsec:othercases},
for the manifold-with-boundary $Z_{K}$. Here is the result:

\begin{proposition}\label{prop:relative-calcuation}
        When $N$ is odd, the integer invariant $q^{w}(Z_{K};\alpha_{0})$ for
        $Z_{K} = S^{1}\times M$ can be expressed in terms of the
    Alexander polynomial $\Delta(t)$ of the knot $K$ by the formula
    \[
            q^{w}(Z_{K};\alpha_{0}) = \pm \prod_{k=1}^{N-1} \Delta({e^{2\pi
                i k/N}}) , 
    \]
     provided this quantity is non-zero. For $k$ non-zero the integer
     invariant $q^{w}(Z_{K};\alpha_{k})$ is zero, because there is no
     moduli space $M^{w}(Z_{K}^{+};\tilde\alpha_{k})$ of formal
     dimension $0$.
\end{proposition}

\begin{proof}

We prove the proposition in a series of lemmas. The main idea is that
we can calculate this invariant by counting certain flat connections
on $Z_{K}$, so reducing the problem to one about representations of
the fundamental group, rather than having to understand non-flat
solutions to the anti-self-duality equations.

\begin{lemma}
    There is a lift $\tilde\alpha_{0} \in \tilde{\Rep}^{w}(T^{3})$ of
    $\alpha_{0} \in \Rep^{w}(T^{3})$ such
    that the moduli space $M^{w}(Z_{K}^{+};\tilde\alpha_{0})$ on the
    cylindrical-end manifold $Z_{K}^{+}$ consists only of flat
    connections and has formal dimension $0$.
\end{lemma}

\begin{proof}
    The moduli space $M^{w}(Z_{K}^{+};\tilde\alpha_{0})$ will consist
    entirely of flat connections if
    $\kappa(Z_{K}^{+};\tilde\alpha_{0})$ is zero, which in turn holds
    if there is at least one flat connection in the moduli space. The
    $\PSU(N)$ representation $\bar\rho_{0}$ of $\pi_{1}(\partial
    Z_{K}; z_{0})$ which corresponds to $\alpha_{0}$ sends $l$ to $1$,
    so it extends as an abelian representation of the fundamental
    group of $Z_{K}$. This representation defines a $\PSU(N)$
    connection on the non-trivial bundle $\adP$, so it defines a lift
    $\tilde\alpha_{0}$ of $\alpha_{0}$ of the sort required.

    By excision, the formal dimension of the moduli space
    $M^{w}(Z_{K}^{+};\tilde\alpha_{0})$ with $\kappa=0$ is independent
    of the knot $K$, so we may consider just the unknot, in which case
    the formal dimension is easily seen to be zero (for example by
    using the product formula for the integer invariants, as in the
    section below).
\end{proof}

\begin{lemma}\label{lem:numberDet}
    The number of flat connections in $M^{w}(Z_{K}^{+};\tilde\alpha_{0})$
    can be expressed in terms of the Alexander polynomial $\Delta(t)$
    by the formula
    \begin{equation}\label{eq:TheAnswer}
                \left|   \prod_{k=1}^{N-1} \Delta({e^{2\pi
                i k/N}}) \right|,
    \end{equation}
    provided that this quantity his non-zero.
\end{lemma}

\begin{proof}
    Because of the previous lemma, we can regard
    $M^{w}(Z_{K}^{+};\tilde\alpha_{0})$ as a moduli space of flat
    connections on the compact space $Z_{K}$. More exactly, it is the
    space of flat $\PSU(N)$ connections in $\adP \to Z_{K}$, divided by the
    determinant-1 gauge group. (Any such flat connection must define
    the representation $\alpha_{0}$ at the boundary, rather than
    $\alpha_{k}$ for non-zero $k$, because the latter have non-zero
    Chern-Simons invariant.)

    We proceed as we did for the $3$-torus $T^{3}$ previously. Let
    $l'$ be a parallel copy of $l$, just as $z'$ was a parallel copy
    of $z$ in our earlier language. We take $l'$ to be the boundary of
    $\Sigma'$, which is a parallel copy of $\Sigma$ in the
    $4$-manifold $Z_{K} = S^{1}\times M$. The surface $\Sigma'$ should
    lie in $p\times M$, where $p$ is distinct from the basepoint $1\in
    S^{1}$.  We equip $w\to Z_{K}$ with a connection $\theta$ that is
    trivial outside the tubular neighborhood of $\Sigma'$ and whose
    curvature is $(-2\pi i )\eta$, where $\eta$ represents the
    Poincar\'e dual class to $\Sigma'$ and is supported in the tubular
    neighborhood.
   Using $\theta$, we can lift a flat connection $A$ in $\adP$ lifts
    to a $U(N)$ connection $\tilde A$ in $P$ that is flat outside the
    neighborhood of $\Sigma'$. The trivialization of $w$ provided by
    $\theta$ reduces the structure group of $P$ to $\SU(N)$ away from
    the tubular neighborhood, and $\tilde A$ becomes a flat $\SU(N)$
    connection there.

    Suppose, then that $A$ and $\tilde A$ are such connections. Let
    \[
                \rho : \pi_{1}(Z_{K}\setminus \Sigma' ; z_{0}) \to
                \SU(N)
    \]
    be the representation defined by the holonomy of $\tilde A$, and
    let
    \begin{equation}\label{eq:barrho}
              \bar \rho : \pi_{1}(Z_{K} ; z_{0}) \to
                \PSU(N)
    \end{equation}
    be defined by the holonomy of $A$. We have already seen, from our
    discussion of the boundary $T^{3}$, that after change of basis we
    must have
    \begin{equation}\label{eq:twoMatrices}
           \begin{aligned}
            \rho(s) &= \begin{bmatrix}
                1 & 0 & 0 & \cdots& 0\\
                0 & \zeta  & 0 & \cdots&0\\
                0 & 0 & \zeta^{2} & \cdots &0\\
                  &   &        &   \ddots&0 \\
               0   &   0&   0 &  0 & \zeta^{N-1}
            \end{bmatrix}
            &
            \rho(m) &= \begin{bmatrix}
                0 & 0 & 0 & \cdots&  \pm 1\\
                1 & 0  & 0 & \cdots&0\\
                0 & 1 & 0 & \cdots &0\\
                  &   &        &   \ddots&0 \\
               0   &   0&   0 &  1 & 0
            \end{bmatrix},
           \end{aligned}
    \end{equation}
    and
    \[
                \rho(l) = \mathbf{1}_{N}.
    \]
    In $\pi_{1}(Z_{K}; z_{0})$
    the class of the circle $s$ is central. So the image of
    $\bar{\rho}$ lies in the centralizer of $\bar{\rho(s)}$. Let
    $\bar{J}$ denote this centralizer, and let $J$ be its inverse
    image in $\SU(N)$. We can describe the group $J$ as a semi-direct product
    \[
            1 \to H \to J \to V \to 1,
    \]
    where $H$ is the standard maximal torus of $\SU(N)$ and $V$ is the
    cyclic group of order $N$ generated by the matrix $\rho(m)$ above.

    We identify the $3$-manifold $M$ with the submanifold $1\times M$
    in the product $Z_{K}  = S^{1}\times M$. Then, by restriction,
    $\rho$ determines a homomorphism
    \[
            \sigma : \pi_{1}(M, z_{0}) \to J,
    \]
    with
    \begin{equation}\label{eq:m-element}
                \sigma(m) = \begin{bmatrix}
                0 & 0 & 0 & \cdots&  \pm 1\\
                1 & 0  & 0 & \cdots&0\\
                0 & 1 & 0 & \cdots &0\\
                  &   &        &   \ddots&0 \\
               0   &   0&   0 &  1 & 0
               \end{bmatrix}
    \end{equation}
    and $\sigma(l) = \mathbf{1}_{N}$.
    Conversely, $\sigma$ determines $\rho$ and hence $\bar\rho$. We
    are therefore left with the task of enumerating such homomorphisms
    $\sigma$.

    Let $\pi$ stand as an abbreviation for the group
    $\pi_{1}(M,z_{0})$, let $\pi^{(1)}$ be the commutator subgroup,
    and $\pi^{(2)}= [\pi^{(1)},\pi^{(1)}]$. The group $J$ is a two-step
    solvable group, with the torus $H$ as commutator subgroup. So we must
    have
    \[
            \sigma(\pi^{(1)}) \subset H,
    \]
    and $\sigma(\pi^{(2)}) = \{1\}$.

    Because $M$ is a knot complement, the quotient $\pi/\pi^{(1)}$ is
    infinite cyclic, generated by the coset represented by the
    meridian $m$.  The abelian group $\pi^{(1)}/\pi^{(2)}$ becomes a
    $\Z[t,t^{-1}]$-module when we let the action of $t$ on
    $\pi^{(1)}/\pi^{(2)}$
    be defined by the action of conjugation, $x \mapsto mxm^{-1}$, on
    $\pi^{(1)}$. In a similar way, the abelian group $H$ becomes a
    $\Z[t,t^{-1}]$-module if we define the action of $t$ to be given by
    conjugation by the element $v \in \SU(N)$ that appears in
    \eqref{eq:m-element}. In this way, $\sigma$ determines (and is
    determined by) a homomorphism of $\Z[t,t^{-1}]$-modules,
    \[
                \sigma' : \frac{\pi^{(1)}}{\pi^{(2)}} \to H.
    \]

    We are left with the task of enumerating such $\Z[t,t^{-1}]$-module
    homomorphisms. The structure of $\pi^{(1)}/\pi^{(2)}$ as a $\Z[t,t^{-1}]$
    module is described by the Alexander polynomial $\Delta$: there is
    an isomorphism
    \begin{equation}\label{eq:AlexanderIdeal}
            \frac{\pi^{(1)}}{\pi^{(2)}} \cong \Z[t,t^{-1}] \bigm/ \Delta(t).
    \end{equation}
   (See \cite{Rolfsen} for example.) So the homomorphism $\sigma'$ is
   entirely determined by the element
   \[
            h =     \sigma'(1) \in H.
   \]
   We switch to additive notation
   for the torus $H$. We write $\tau_{v}$ for the automorphism of $H$
   given by conjugation by $v$. In the ring of endomorphisms of $H$,
   we can consider the element $\Delta(\tau_{v})$. The description
   \eqref{eq:AlexanderIdeal} of $\pi^{(1)}/\pi^{(2)}$ means that
   the element $h = \sigma'(1)$ must lie in the kernel of
   $\Delta(\tau_{v})$.

   Thus we have shown that the moduli space
   $M^{w}(Z_{K}^{+};\tilde\alpha_{0})$ is in one-to-one correspondence
   with the elements $h$ in the abelian group
   \[
        \ker(\Delta(\tau_{v})) \subset H.
   \]
   Let $\tilde H$ be the universal cover of $H$ (or the Lie algebra of
   this torus), and let $\tilde\tau_{v}$ be the lift of $\tau_{v}$ to
   this vector space. Then the order of the kernel of
   $\Delta(\tau_{v})$ is equal to the absolute value of the
   determinant
   \[
            \det \bigl ( \Delta(\tilde \tau_{v})\bigr) 
   \]
   of the linear transformation $\Delta(\tilde \tau_{v})$, provided
   that this determinant is non-zero. If the determinant is zero, then
   the kernel is infinite, for it is a union of tori, whose dimension
   is equal to the dimension of the null-space of $\Delta(\tilde
   \tau_{v})$. As a linear transformation of the $(N-1)$-dimensional
   real vector space $\tilde H$, the operator $\tilde
   \tau_{v}$ is diagonalizable after complexification, and its
   eigenvalues on $\tilde H \otimes \C$ are the non-trivial $N$-th
   roots of unity, $\zeta^{k}$ for $k=1,\dots, N-1$. The eigenvalues
   of $\Delta(\tilde
   \tau_{v})$ are therefore the complex numbers $\Delta(\zeta^{k})$,
   so
   \[
        \det \bigl ( \Delta(\tilde \tau_{v})\bigr) =
        \prod_{k=1}^{N-1}\Delta(\zeta^{k}).
   \]
   This completes the proof of the lemma.
\end{proof}

\begin{remark}
    The formula \eqref{eq:TheAnswer} has another, closely-related
    interpretation. It is the order of the first homology group
    $H_{1}(Y_{N}; \Z)$ of the $3$-manifold $Y_{N}$ obtained as the
    $N$-fold cyclic cover of $S^{3}$ branched over the knot. By
    Poincar\'e duality, this group is the same as $H^{2}(Y_{N};\Z)$,
    which (when finite) classifies flat complex line bundles on
    $Y_{N}$. Let $p:M'\to M$ be the $N$-fold cyclic cover of the knot
    complement. If we take a flat complex line bundle on $Y_{N}$ and
    restrict it to $M'$, we can then
    push it forward by the map $p$ to obtain a flat rank
    $N$ bundle on $M$. In the case that $N$ is odd, this construction
    produces a
    flat $\SU(N)$ bundle on $M$ for each element of $H^{2}(Y_{N};\Z)$.
    In another language, we are constructing $N$-dimensional
    representations of $\pi_{1}(M$) as induced representations,
    starting from $1$-dimensional representations of the subgroup
    $\pi_{1}(M')$.
    (In the case $N$ even, a slight adjustment is needed,
    corresponding to the negative sign in the top-right entry of
    $\rho(m)$ above.)
\end{remark}

\begin{lemma}\label{lem:regularFlat}
    If the quantity in Lemma~\ref{lem:numberDet} is non-zero, then the
    moduli space $M^{w}(Z_{K}^{+};\tilde\alpha_{0})$ is regular.
\end{lemma}

\begin{proof}
    Let $A$ be a flat connection in $\adP$ representing a point of
    this moduli space. The formal dimension of the moduli space is
    zero, and $A$ is irreducible. So to establish regularity, it is
    sufficient to show that $H^{1}_{A}$ is zero, or equivalently that
    $\ker d^{+}_{A}$ is equal to $\mathrm{im} d_{A}$ in
    $L^{2}_{l,A}(Z^{+}_{K}; \g_{P})$. Harmonic theory on the
    cylindrical end manifold means that we can rephrase this as saying
    that the kernel of $d^{+}_{A} \oplus d^{*}_{A}$ is zero. Any
    element of the kernel of this Fredholm operator is exponentially
    decaying on the cylindrical end, and we can integrate by parts
    (using the fact that $A$ is flat), to conclude that an element $a$
    in the kernel must have $d^{-}_{A}a=0$.  Thus $a$ is in the kernel
    of the operator $d_{A}\oplus d^{*}_{A}$.

    The kernel of $d_{A}\oplus d^{*}_{A}$ represents the first
    cohomology group with coefficients in the flat bundle $\g_{P}$
    with connection $A$:
    \[
                H^{1}(Z_{K}; (\g_{P},A)).
    \]
    There is no difference here between the absolute or relative
    group, because the connection $A$ is irreducible on the boundary
    of $Z_{K}$. We must therefore only show that this cohomology group
    is zero.

    Under the adjoint action of the matrix $\rho(s)$ from
    \eqref{eq:twoMatrices}, the complexified Lie algebra $\sl(N,\C)$
    decomposes as a sum
    \[
                \sl(N,\C) = E_{0} \oplus \cdots \oplus E_{N-1},
    \]
    where $E_{k}$ is the $\zeta^{k}$-eigenspace of
    $\mathrm{Ad}(\rho(s))$, which is the span of the elementary
    matrices $e_{ij}$ with a $1$ in row $i$ and column $j$, with $i-j
    = k \pmod{N}$.  There is a corresponding decomposition of the
    bundle $\g_{P}$ which is respected by the flat connection, because
    $s$ is central in the fundamental group.

    So we look now at the individual summands $H^{1}(Z_{K}; (E_{k},
    A))$. The manifold $Z_{K}$ is a product, and $s$ is the element of
    $\pi_{1}$ corresponding to the circle fiber. On the circle
    $H^{1}(S^{1}, \lambda)$ is zero if $\lambda$ is a local system
    with fiber $\C$ and non-trivial holonomy. By the Kunneth theorem,
    it follows that $H^{1}(Z_{K}; (E_{k},
    A))$ is zero for $k$ non-zero.

    We are left with the summand $H^{1}(Z_{K}; (E_{0},
    A))$. The holonomy of $A$ on $E_{0}$ is trivial in the circle
    factor of $Z_{K}$; so we may as well consider $H^{1}(M;
    (E_{0},A))$. On the knot complement $M$, the bundle $E_{0}$ has
    rank $N-1$ and the holonomy of $A$ factors through the
    abelianization of the fundamental group, $H_{1}(M;\Z)$, which is
    generated by $m$. Under the action of $\mathrm{Ad}(\rho(m))$, the
    bundle $E_{0}$ decomposes into $N-1$ flat line bundles $L_{k}$
    ($k=1, \dots N$); and the holonomy of $L_{k}$ along the meridian
    is $\zeta^{k}$.

    To complete the proof, we have to see that $H^{1}(M; (L_{k},A))$
    is zero unless $\Delta(\zeta^{k})=0$. This is a standard result,
    and can be proved as follows. If the cohomology group is non-zero,
    let $\alpha$ be a non-zero element, and let $\alpha'$ be the
    pull-back to the $N$-fold cyclic cover $p: M' \to M$. This
    pull-back is an element of $H^{1}(M'; \C)$ and has trivial
    pairing with the loop $m' = p^{-1}(m)$. It therefore pairs
    non-trivially with a loop $\gamma\subset M'$ that lifts to the
    $\Z$-covering $\tilde M$. Thus the pull-back of $\alpha$ to the
    $\Z$-covering is a non-zero class $\tilde\alpha\in H^{1}(\tilde
    M; \C)$. By construction $\tilde\alpha$ lies in the
    $\zeta^{k}$-eigenspace of the action of the covering
    transformation. Up to units in the ring $\C[t,t^{-1}]$, the
    characteristic polynomial of the covering transformation acting on
    $H^{1}(\tilde M; \C)$ is $\Delta(t)$. So $\Delta(\zeta^{k})$ is
    zero.
\end{proof}

\begin{lemma}\label{lem:same-sign}
    Suppose $N$ is odd. Then under the same hypothesis as
    Lemma~\ref{lem:regularFlat}, all the points of the moduli space
    $M^{w}(Z_{K}^{+};\tilde\alpha_{0})$ have the same sign.
\end{lemma}

\begin{proof}
We wish to compare the orientations of two points $[A_{0}]$ and
$[A_{1}]$ in the moduli space $M^{w}(Z_{K}^{+};\tilde\alpha_{0})$. Let
$[A_{t}]$ be a 1-parameter family of connections in
$\cA(Z^{+}_{K};\tilde\alpha_{0})$ joining $A_{0}$ to $A_{1}$. Let
$D_{t}$ be the corresponding Fredholm operators on the cylindrical-end
manifold $Z^{+}_{K}$, as in \eqref{eq:Fredholm-cylinder}. The
operators $D_{0}$ and $D_{1}$ are invertible, and so there are
canonical trivializations of the determinant lines $\det(D_{0})$ and
$\det(D_{1})$. The lemma asserts that these canonical trivializations
can be extended to a trivialization of the determinant line
$\det(D_{\bullet})$ on $[0,1]$ .

Let $p: M' \to M$ be again the $N$-fold cyclic cover of the knot
complement, and let $p : Z'_{K} \to Z_{K}$ be the corresponding cyclic
cover of $Z_{K} = S^{1}\times M$. Let ${Z'}^{+}_{K}$ be the
corresponding cylindrical-end manifold. Set
\[
            A'_{t} = p^{*}(A_{t})
\]
These are connections on ${Z'}^{+}_{K}$, asymptotic
to $p^{*}(\alpha_{0})$ on the cylindrical end. There are corresponding
Fredholm operators $\det(D'_{t})$, for $t\in [0,1]$; we use a weighted
Sobolev space with a small exponential weight to define the domain and
codomain of $D'_{t}$, because the pull-back of $\tilde\alpha_{0}$ under
the covering map $p$ is not isolated in the representation variety.
The operators $D'_{0}$ and $D'_{1}$ are invertible, as one can see by
a calculation similar to that in the proof of
Lemma~\ref{lem:regularFlat} above.

The covering group $\Z/(N\Z)$ acts on the domain and codomain of
$D'_{t}$, which therefore decompose into isotypical components
according to the representations of this cyclic group. The restriction of
$D'_{t}$ to the component on which $\Z/(N\Z)$ acts trivially is just
$D_{t}$ (with the inconsequential change that $D_{t}$ is acting now
on a weighted Sobolev space). The non-trivial representations of
$\Z/(N\Z)$ are all complex, because $N$ is odd, and their contribution
to the determinant is therefore trivial. We can therefore replace
$D_{t}$ by $D'_{t}$ without changing the question: we want to know
whether the canonical trivializations of $\det(D'_{0})$ and
$\det(D'_{1})$ can be extended to a trivialization of
$\det(D'_{\bullet})$ on the interval $[0,1]$.

Let $\tilde A'_{t}$ be the $U(N)$ connection with determinant
$p^{*}(\theta)$, corresponding to the $\PSU(N)$ connection $A'_{t}$.
The description of the representations $\bar\rho$ corresponding to the
flat connections in the moduli space
$M^{w}(Z^{+}_{K};\tilde\alpha_{0})$ show that $\tilde A'_{0}$ and
$\tilde A'_{1}$ are each
compatible with a decomposition of the $U(N)$ bundle $P' = p^{*}(P)$
into a direct sum of $U(1)$ bundles (though a different decomposition
in the two cases). That is, we have
decompositions of the associated rank-$N$ vector bundle $E$ as
\[
            E = L_{i,1} \oplus \dots \oplus L_{i,N}
\]
for $i=0$ and $1$. Since $\tilde A'_{i}$ is projectively flat, the
first Chern class of each $L_{i,m}$ is the same over the reals,
independent of $i$ and $m$: the common first Chern class of all these
line bundles is $(1/N)c_{1}(p^{*}(w))$.  We can include $Z_{K}'$ in a
closed $4$-manifold $W$, and since all the $L_{i,m}$ are the same on
the boundary, we can suppose they all extend to $W$ in such a way that
their rational first Chern classes are equal. The comparison of the
orientations determined by $D'_{0}$ and $D'_{1}$ is then the same as
the comparison between orientations
\[
            O(L_{i,1}, \dots, L_{i,n}, o_{W}) , \qquad (i=0,1)
\]
on $W$. The analysis of section~\ref{subsec:OrientationConventions}
shows that these two orientations depend only on the real or rational
first Chern classes of the line bundles, so the two orientations
agree.
\end{proof}

The proof of Proposition~\ref{prop:relative-calcuation} is now
complete, because the definition of the invariant is the signed count
of the points of the moduli space.
\end{proof}

\begin{remarks}
    When $N$ is even, the proof of Lemma~\ref{lem:same-sign} breaks
    down: in the decomposition of the domain and codomain of $D'_{t}$,
    there is now a component on which the generator of  $\Z/N\Z$ acts
    as $-1$; and we know no more about the determinant line of $D'_{t}$ on
    this component than we did about the determinant line of the
    original $D_{t}$.

    One should expect to prove
    Proposition~\ref{prop:relative-calcuation} also in the case that
    the expression that appears there is zero. In this case, the
    moduli space is a union of circles or higher-dimensional tori, and
    one should try to use a holonomy perturbation to make the perturbed
    moduli space empty.
\end{remarks}

\subsection{The Fintushel-Stern construction}

Let $X$ be a closed oriented $4$-manifold with $b^{+}_{2}(X)\ge 2$. Let $w$ be
a line bundle with $c_{1}(w)$ coprime to $N$, and suppose that the integer
invariant $q^{w}(X)$ is non-zero.

Let $T$ be an embedded torus in $X$, with trivial normal bundle, and
suppose that the pairing of
$c_{1}(w)$ with $[T]$ is $1$. The torus $T$ has a closed tubular neighborhood $N\subset
X$ of the form $T\times D^{2}$. If $U$ is the unknot in $S^{3}$, then
the $4$-manifold $Z_{U}$ obtained from the unknot is also $T\times
D^{2}$. We choose an identification $\phi$ of $Z_{U}$ with $N$, in such a way
that the longitudinal curve $l$ on the boundary of $Z_{U}$ is matched
by $\phi$
with the curve $\mathrm{(point)} \times \partial D^{2}$ that links the
embedded torus. Using the standard curves $l$, $m$ and $s$, the
boundary of $Z_{U}$ is canonically identified with the boundary of
$Z_{K}$, for any other knot $K$ in $S^{3}$. Thus the choice of $\phi$
gives us preferred diffeomorphisms
\[
            \psi : \partial Z_{K} \to \partial N
\]
for all $K$. We now form a closed $4$-manifold $X_{K}$ by removing $N
= Z_{U}$ and replacing it with $Z_{K}$:
\[
            X_{K} = Z_{K} \cup_{\psi} (X\setminus N).
\]
This is Fintushel and Stern's construction, from
\cite{Fintushel-Stern-Knot}. The Mayer-Vietoris sequence tells us that
the cohomology rings of $X$ and $X_{K}$ are canonically isomorphic. 

The line bundle $w$ on $X\setminus N$
extends to $X_{K}$. The extension $w$ is unique if we demand the  additional property that
$c_{1}(w)^{2}[X_{K}]$ and $c_{1}(w)^{2}[X]$ are equal.
The moduli spaces $M^{w}_{\kappa}(X_{K})$ and $M^{w}_{\kappa}(X_{K})$
have the same formal dimension, and there is a potentially non-zero
integer invariant $q^{w}(X_{K})$.  We can calculate this invariant for
odd $N$:

\begin{proposition}
    When $N$ is odd, the integer invariant $q^{w}(X_{K})$ can be expressed in terms of
    $q^{w}(X)$ and the Alexander polynomial $\Delta$ of $K$ as
    \[
            q^{w}(X_{K}) =  q^{w}(X) \times
            \prod_{k=1}^{N-1}\Delta(\zeta^{k}),
    \]
    where $\zeta = e^{2\pi i /N}$.
\end{proposition}

\begin{proof}
    Up to an overall ambiguity in the sign, this proposition
    is now a formal consequence of the calculation of the
    invariants $q^{w}(Z_{K};\alpha_{k})$ and $q^{w}(Z_{U};\alpha_{k})$
    which are provided by Proposition~\ref{prop:relative-calcuation},
    and the product law \eqref{eq:generalProduct}. These tell us that
    \[
             q^{w}(X_{K})  = \pm\  n \times
             \prod_{k=1}^{N-1}\Delta(\zeta^{k}),
    \]
    where $n = q^{w}(X\setminus N^{\circ};\alpha_{0})$ is a quantity
    independent of $K$. Because $X_{U}$ is $X$, and the Alexander
    polynomial of the unknot is $1$, it follows that $n = \pm
    q^{w}(X)$.

    It remains to check the overall sign. Certainly, the sign is
    correct if $K$ is the unknot. Note also that since the complex
    numbers $\Delta(\zeta^{k})$ come in conjugate pairs when $N$ is
    odd, their product is positive; so our task is to check that the
    sign of the invariant is independent of $K$.  To this end, we
    observe that among the flat connections on the cylindrical-end
    manifold $Z_{K}^{+}$ that are enumerated in
    Lemma~\ref{lem:numberDet}, there is a unique connection, say
    $[A_{*}]$, distinguished by the fact that the corresponding
    representation $\bar\rho$ of the fundamental group (see
    \eqref{eq:barrho}) factors through the abelianization
    $H_{1}(Z_{K};\Z)$. Let $[B]$ be an isolated solution in the moduli
    space $M^{w}_{\bomega}((X\setminus
    N^{\circ})^{+};\tilde\alpha_{0})$, and let $[C]$ be the
    connection
    on $X_{K}$ obtained as the result of gluing $[B]$ to
    $[A_{*}]$ with a long neck.  It is not important that $C$ is
    actually a solution of the equations: we can instead just use
    cut-off functions to patch together the connections $A_{*}$ and
    $B$; but it is important that the operator $D_{C}$ is invertible
    \eqref{eq:D-op},
    so that we can ask about the sign of the point $[C]$.

    To complete the proof we will show that the sign of $[C]$ is
    independent of the knot $K$. If $K$ is any knot and $K_{0}$ is the
    unknot, then we can transform
    $Z_{K} = S^{1}\times M$ into $Z_{K_{0}}$ by a sequence of surgeries along tori
    \[
         T =   S^{1}\times\delta \subset S^{1} \times M,
    \]
    where in each surgery the curve $\delta$ is a null-homologous
    circle in $M$ (compare \cite{Fintushel-Stern-Knot}). The
    connection $A_{*}$ (and therefore also $C$) is flat on such a
    torus and lifts to an abelian representation $\rho$ of
    $\pi_{1}(T)$, sending the $S^{1}$ generator to the first matrix in
    \eqref{eq:twoMatrices} and the generator $\delta$ to $1$.  Thus
    $C$ is compatible with a decomposition of the associated rank-$N$
    vector bundle as a sum of flat line bundles. 
    It now follows from excision and the material of
    section~\ref{subsec:OrientationConventions} that the surgery does
    not alter the sign of $[C]$.
\end{proof}

As a particular case, we can take $X$ to be a $K3$ surface. With $w$
as in section~\ref{sec:K3}, we calculated the invariant $q^{w}(X)$ to
be $1$, for the preferred homology-orientation of $X$. The value of
$c_{1}(w)$ in our calculation was $(N+1)h$; but we could equally have
taken $c_{1}(w)$ to be just the primitive class $h$ (with square
$2(N-1)$). In this case, the invariant $q^{w}(X)$ would have again
been $1$, by Proposition~\ref{prop:liftcompare}.  A primitive
cohomology class on a $K3$ surface has pairing $1$ with some embedded
torus $T$. So choose $T$ with $h[T]=1$. We are then in a position to
apply the proposition above, to obtain
\[
q^{w}(X_{K})  = 
             \prod_{k=1}^{N-1}\Delta(\zeta^{k}).
\]             
for odd $N$. As Fintushel and Stern observe in \cite{Fintushel-Stern-Knot}, the
$4$-manifold $X_{K}$ is a homotopy $K3$-surface if we take $T$ to be a
standardly embedded torus (so that $X\setminus T$ is simply
connected).

\section{Polynomial invariants}
\label{sec:invariants}

We now generalize the integer-valued invariant $q^{w}(X)$, to define
polynomial invariants, as Donaldson did in \cite{Donaldson-polynomials}. Our
approach follows \cite{Donaldson-polynomials},
\cite{Donaldson-Kronheimer} and \cite{KM-Structure} closely.

\subsection{Irreducibility on open sets}

Let $\Omega\subset X$ be a connected, non-empty open set, chosen so that
$\pi_{1}(\Omega;x_{0}) \to \pi_{1}(X;x_{0})$ is surjective, for some
$x_{0}\in \Omega$. The unique continuation argument from
\cite{Donaldson-Kronheimer} shows that if an anti-self-dual connection
$[A] \in M^{w}_{\kappa}(X)$ is irreducible, then its restriction to
$\Omega$ is also irreducible. For the perturbed equations, it may be
that this result fails. 
The next lemma provides a suitable substitute for our purposes. It is
only a very slight modification of Proposition~\ref{prop:smallPert}, and we
therefore omit the proof.

\begin{lemma}\label{lem:smallPert-Omega}
    Let $\Omega\subset X$ be as above.
    Let $\kappa_{0}$ be given, and let $g$ be a metric on $X$ with the property that the unperturbed
    moduli spaces $M^{w}_{\kappa}$ contain no reducible solutions for
    any $\kappa \le \kappa_{0}$. Then there exists $\epsilon > 0$ such
    that for all $\bomega\in W$ with $\| \bomega \|_{W} \le \epsilon$,
    and all $[A]$ in $M^{w}_{\kappa,\bomega}(X)$, the restriction of
    $A$ to $\Omega$ is irreducible. \qed
\end{lemma}

\subsection{Cohomology classes and submanifolds}

Fix as usual a $U(N)$ bundle $P\to X$, and let $\B^{*}=\B^{*}(X)$ denote the
irreducible part of the usual configuration space. There is a
universal family of connections, carried by a bundle
\[
            \ad(\mathbb{P}) \to X\times \B^{*}.
\]
This universal bundle can be constructed as the quotient of $\adP
\times \cA^{*}$ by the action of $\G$. Note, however, that we cannot
form a universal $U(N)$ bundle $\mathbb{P}$ in this way, and our
notation is not meant to imply that such a $U(N)$ bundle exists. We
define a $4$-dimensional cohomology class on $X\times \B^{*}$ by taking the
first Pontryagin class of the adjoint bundle, with the now-familiar
normalization:
\[
\begin{aligned}
    \cc &= -(1/2N) p_{1}( \su_{\mathbb{P}} ) \\
       &=  (1/2N) c_{2}( \sl_{\mathbb{P}} ) \\
       &\in H^{4}( X\times\B^{*};\Q).
\end{aligned}
\]
Using the slant product, we now define
\begin{equation}\label{eq:mu-map}
        \mu: H_{i}(X;\Q) \to H^{4-i}(\B^{*};\Q)
\end{equation}
by the formula $\mu(\alpha) = \cc / \alpha$. We will concern ourselves
here only with even-dimensional classes on $X$,
so we define $\bbA(X)$ to be the polynomial algebra
\[
    \bbA(X)  = \mathrm{Sym}( H_{\mathrm{even}}(X;\Q) ),
\]
and we extend $\mu$ to a ring homomorphism
\[
        \mu: \bbA(X) \to H^{*}(\B^{*};\Q).
\]
We regard $\bbA(X)$ as a graded algebra, defining the grading of
$H_{i}(X;\Q)$ to be $4-i$ so that $\mu$ respects the grading.

Now let $\Sigma$ be an oriented embedded surface in $X$ representing a
2-dimensional homology class $[\Sigma]$. Although $\mu([\Sigma])$ is
not necessarily an integral cohomology class, the definition of $\mu$
makes evident that $(2N)\mu([\Sigma])$ is integral. There is therefore
a line bundle
\[
            \cL \to B^{*}
\]
with $c_{1}(\cL) = (2N)\mu([\Sigma])$. We can realize this line bundle
as a determinant line bundle: if $\bar\partial$ denotes the
$\bar\partial$-operator on $\Sigma$, then
\[
            \cL = \det ( \bar\partial^{*}_{\su(\mathbb{P})} ),
\]
where the notation means that we couple $\bar\partial^{*}$ to the
$\PSU(N)$ connections in the family over $\B^{*}$. The following lemma
and Corollary~\ref{cor:transverseVsigma} below play the same role as
Lemma~(5.2.9) of \cite{Donaldson-Kronheimer}. There is a slight extra
complication in our present setup: we need to deal with the weaker
$L^{p}_{1}$ norms, because of the way our perturbations behave under
Uhlenbeck limits. (See Proposition~\ref{prop:UhlenbeckPert}.)

\begin{lemma}\label{lem:sections}
    Fix an even integer $p > 2$. Then
    there is a complex Banach space $E$ and a continuous linear map $s : E \to
    \Gamma(\B^{*};\cL)$ with the following properties.
    \begin{enumerate}
       \item For each $e\in E$, the section $s(e)$ is smooth, and
        furthermore the map \[
         s^{\dag}: E \times \B^{*}\to \cL \]
        obtained by evaluating $s$ is a smooth map of Banach
        manifolds.

        \item For each $[A]\in \B^{*}$, the exists $e\in E$ such that
        $s(e)([A])$ is non-zero.

        \item For each $e\in E$, the section $s(e)$ is continuous in
        the $L^{p}_{1}$ topology on $\B^{*}$.
    \end{enumerate}
\end{lemma}

\begin{proof}
    Let $[A]$ belong to $\B^{*}$, and let $S \subset \cA$ be the
    Coulomb slice through $A$. Using the definition of $\cL$ as a
    determinant line bundle, construct a non-vanishing local section
    $t$ for $\cL$ on some neighborhood $U$ of $A$ in $S$. The section
    $t$ can be constructed to be smooth in the $L^{q}_{1}$ topology
    for all $q$; so we can take it that $U$ is open in the $L^{p}_{1}$
    topology. Because $p$ is an even integer, the $L^{p}_{1}$ norm
    (or an equivalent norm) has the property that $\| \;
    \|^{p}_{L^{p}_{1}}$ is a smooth function. We can therefore
    construct a smooth cut-off $\beta$ function on the $L^{p}_{1}$ completion of $S$,
    supported in (the completion) of $U$. Then the section $s = \beta
    t$ can be extended as a smooth section of $\cL$ on the $L^{p}_{1}$
    configuration space $\B^{*}_{L^{p}_{1}}$.

    Now choose a countable section collection of points $[A_{i}]$ in
    $\B^{*}$ and sections $s_{i}$ constructed as above, such that the
    support of the $s_{i}$ covers $\B^{*}$.  Fix a collection of
    positive real numbers $c_{i}$, and let $E$ be the Banach space of all
    sequences $\mathbf{z}=\{z_{i}\}_{i\in \mathbb{N}}$ such that the sum
    \[
                \sum c_{i} |z_{i}|
    \]
    is finite. If the sequence $c_{i}$ is sufficiently rapidly
    increasing, then for all $\mathbf{z}$ in $E$, the sum
    \[
            \sum z_{i} s_{i}
    \]
    converges to a smooth section of $\cL$ on $\B^{*}$, and the
    resulting linear map $s : E \to \Gamma(\B^{*};\cL)$ satisfies the
    conditions of the lemma.
\end{proof}

\begin{corollary}\label{cor:transverseVsigma}
    Let $M$ be a finite-dimensional manifold and $f : M \to \B^{*}$ a
    smooth map. Then there is a section $s$ of $\cL \to \B^{*}$ such
    that the section $s\circ f$ of the line bundle $f^{*}(\cL)$ on $M$
    is transverse to zero. \qed
\end{corollary}

The above corollary can be applied directly to the situation that $M$
is a regular moduli space $M^{w}_{\kappa,\bomega}$ of perturbed
anti-self-dual connections, containing no reducibles. In this case,
the zero set of $s$ on $M$ is a cooriented submanifold of codimension
$2$ whose dual class is $(2N) \mu([\Sigma])$.  As with Donaldson's
original construction however, we must modify the construction to
obtain a codimension-2 submanifold which behaves a little better with
respect to the Uhlenbeck compactification.

Let $\Sigma$ be an embedded surface as above, and let $\nu(\Sigma)$ be
a closed neighborhood that is also a 4-dimensional submanifold with
boundary in $X$. By adding $1$-handles to $\nu(\Sigma)$, we may
suppose that this neighborhood is connected and that
$\pi_{1}(\nu(\Sigma)) \to \pi_{1}(X)$ is surjective.  Let
$\B^{*}(\nu(\Sigma)$ be the configuration space of irreducible
connections on $\nu(\Sigma)$. We suppose $w$ and $g$ are such that the
moduli spaces $M^{w}_{\kappa'}$ contain no reducible solutions for
all $\kappa'\le \kappa$. Lemma~\ref{lem:smallPert-Omega} tells us
that we may choose a small regular perturbation $\bomega$ so that
there are well-defined restriction maps
\[
            f : M^{w}_{\bomega,\kappa'} \to \B^{*}(\nu(\Sigma)).
\]
The determinant line bundle $\cL$ is defined on $\B^{*}(\nu(\Sigma))$.
As in \cite{Donaldson-Kronheimer}, we choose a section $s$ of the
$\cL \to \B^{*}(\nu(\Sigma))$ such that the section $s\circ f$ is
transverse to zero. This can be done by repeating the proof of the
lemma above with $\B^{*}(\nu(\Sigma))$ replacing $\B^{*}(X)$. We can
even choose $s$ so that transversality holds for all $\kappa' \le
\kappa$ simultaneously.

As is customary, we write $V_{\Sigma} \subset \B^{*}(\nu(\Sigma))$ for
the zero-set of a section $s$ of this sort; and we write
\[
            M^{w}_{\bomega,\kappa} \cap V_{\Sigma}
\]
for the zero set in the moduli space, suppressing mention of the
restriction map $f$.

\begin{lemma}\label{lem:dichotomy}
    Let $V_{\Sigma}$ be constructed as above, and let
    $[A_{n}]$ be a sequence in $M^{w}_{\bomega,\kappa} \cap V_{\Sigma}$,
    converging in the sense of Proposition~\ref{prop:UhlenbeckPert} to
    an ideal solution $([A'], \bx)$, with $[A']$ in
    $M^{w}_{\bomega,\kappa-m}$. Then either:
    \begin{enumerate}
        \item the limit $[A']$ belongs to $M^{w}_{\bomega,\kappa-m} \cap
        V_{\Sigma}$; or
        \item one of the points of the multi-set $\bx$ lies in the
        closed neighborhood $\nu(\Sigma)$.
    \end{enumerate}
\end{lemma}

\begin{proof}
    If $\bx$ is disjoint from $\nu(\Sigma)$, then the $A_{n}$ converge
    to $A'$
    in the $L^{p}_{1}$ topology on $\nu(\Sigma)$ after gauge
    transformation. The set $V_{\Sigma}$ is closed in
    $\B^{*}(\nu(\Sigma))$ in the $L^{p}_{1}$ topology, because of the
    last condition in Lemma~\ref{lem:sections}, so $[A']$ belongs to
    $V_{\Sigma}$.
\end{proof}

We now wish to represent $(2N)\mu([x])$ ($x$ is a point in $X$) by a
codimension-4 intersection $M^{w}_{\bomega,\kappa} \cap V_{x}$ in the
same manner. As in \cite{KM-Structure}, we must now take $V_{x}$ to be
a stratified space with smooth manifolds for strata, rather than
simply a smooth submanifold.

To obtain a suitable definition of $V_{x}$, we observe that $(2N)
\mu([x])$ is represented on $\B^{*}$ by $c_{2}(\sl_{\mathbb{P_{x}}})$,
where $\mathbb{P}_{x}\to\B^{*}$ is the base-point fibration: the
restriction of $\mathbb{P}$ to $\{x\} \times \B^{*}$. To represent
$c_{2}$ geometrically, we take $N^{2}-2$ sections $s_{1}, \dots
s_{N^{2}-2}$ of the complex vector bundle $\sl_{\mathbb{P}_{x}}$ and
examine the locus on which they are linearly dependent. We
must require a transversality condition: the following lemma and
corollary are essentially the same as Lemma~\ref{lem:sections} and
Corollary~\ref{cor:transverseVsigma} above.

\begin{lemma}\label{lem:multisections}
    Fix an even integer $p > 2$, and write $\mathcal{V} \to \B^{*}$ for the
    bundle $\Hom(\C^{N^{2}-2}, \sl_{\mathbb{P}_{x}})$. Then
    there is a complex Banach space $E$ and a continuous linear map $s : E \to
    \Gamma(\B^{*};\mathcal{V})$ with the following properties.
    \begin{enumerate}
       \item For each $e\in E$, the section $s(e)$ is smooth, and
        furthermore the map \[
         s^{\dag}: E \times \B^{*}\to \cL \]
        obtained by evaluating $s$ is a smooth map of Banach
        manifolds.

        \item For each $[A]\in \B^{*}$, the map $E \to
        \mathcal{V}_{[A]}$ obtained by evaluating $s$ is surjective.

        \item For each $e\in E$, the section $s(e)$ is continuous in
        the $L^{p}_{1}$ topology on $\B^{*}$.
    \end{enumerate} \qed
\end{lemma}

\begin{corollary}\label{cor:transverseVx}
    Let $M$ be a finite-dimensional manifold and $f : M \to \B^{*}$ a
    smooth map. Then there is a section $s$ of\/ $\Hom(\C^{N^{2}-2},
    \sl_{\mathbb{P}_{x}}) \to \B^{*}$ such
    that the section $s\circ f$ of the pull-back $f^{*}(\cL)$ on $M$
    is transverse to the stratification of\/  $\Hom(\C^{N^{2}-2},
    \sl_{\mathbb{P}_{x}})$ by rank. \qed
\end{corollary}

As with the case of $2$-dimensional classes, we now take a closed
neighborhood $\nu(x)$ that is also a manifold with boundary. We require
that $\pi_{1}(\nu(x), x) \to \pi_{1}(X,x)$ is surjective. We can again
suppose that for $\kappa'\le \kappa$, the restriction of every solution in
$M^{w}_{\bomega,\kappa'}$ to $\nu(x)$ is irreducible, so that we have
maps
\[
            f : M^{w}_{\bomega,\kappa'} \to \B^{*}(\nu(x)).
\]
We take a section $s$ satisfying the transversality condition of the
above corollary, taking $M$ to be the union of the
$M^{w}_{\bomega,\kappa'}$ and replacing $\B^{*}$ by $\B^{*}(\nu(x))$
in the statement. Write 
\[ V_{x}\subset \B^{*}(\nu(x))
\]
for the locus where $s$ does not have full rank. Then the intersection
\[
        M^{w}_{\bomega,\kappa}\cap V_{x}
\]
is a stratified subset, with each stratum a smooth, cooriented
manifold, and top stratum of real codimension $4$ in
$M^{w}_{\bomega,\kappa}$. The dual class to $V_{x}$ is $(2N)\mu([x])$.
Lemma~\ref{lem:dichotomy} continues to hold with $V_{x}$ in place of
$V_{\Sigma}$.

\subsection{Polynomials in the coprime case}

We now have the necessary definitions to construct the $\PSU(N)$
generalization of Donaldson's polynomial invariants, in the case that
$w = c_{1}(P)$ is coprime to $N$.

We suppose as usual that $b^{+}_{2}(X)$ is at least $2$. Because of
the coprime condition, we can choose a Riemannian metric $g$ on $X$
such that the moduli spaces $M^{w}_{\kappa}$ contain no reducible
solutions. Given some $\kappa_{0}$, we may find $\epsilon$ so that for
all perturbations $\bomega$ with $\|\bomega\|_{W}\le \epsilon$, the
moduli spaces $M^{w}_{\bomega,\kappa}$ contain no reducibles for
$\kappa \le \kappa_{0}$. Amongst such small perturbations, we may
choose one so that the moduli spaces are regular. If $N$ is even,
choose a homology orientation for $X$ so as to orient the moduli
spaces.

Let $d$ be the dimension of $M^{w}_{\bomega,\kappa}$, and suppose that
$d$ is even. (This parity condition holds if $N$ is odd, or if
$b^{+}_{2}(X) - b_{1}(X)$ is odd.)  Consider an element of $\bbA(X)$,
given by a monomial
\begin{equation}\label{eq:z-monomial}
        z = [x]^{r} [\Sigma_{1}][\Sigma_{2}] \dots [\Sigma_{t}].
\end{equation}
Suppose that $4r + 2t = d$, so that $z$ is of degree $d$. We may
assume that the surfaces $\Sigma_{i}$ intersect only in pairs, and we
can represent the class $[x]$ by $r$ distinct points $x_{1}$, \dots,
$x_{r}$ all disjoint for the surfaces. Choose neighborhoods
\[
            \nu(x_{1}), \dots,  \nu(x_{r}) , \nu(\Sigma_{1}) , \dots ,
            \nu(\Sigma_{t})
\]
satisfying the surjectivity condition on $\pi_{1}$ as above. We may
assume that these satisfy the same intersection conditions: the only
non-empty intersections are possible intersections of pairs
$\nu(\Sigma_{i}) \cap \nu(\Sigma_{j})$. Construct submanifolds (or
stratified subspaces) $V_{x_{i}}$ and $V_{\Sigma_{j}}$ as above. By a
straightforward extension of the transversality argument, we may
arrange that all multiple intersections
\[
            M^{w}_{\bomega,\kappa'} \cap V_{x_{i_{1}}} \cap
            \cdots \cap V_{x_{i_{p}}}\cap V_{\Sigma_{j_{1}}} \cap
            \cdots \cap V_{\Sigma_{j_{q}}}
\]
are transverse. In the case of the $V_{x_{i}}$, transversality means
transversality for each stratum. The intersection
\begin{equation}\label{eq:q-intersection}
            M^{w}_{\bomega,\kappa} \cap V_{x_{1}} \cap
            \cdots \cap V_{x_{r}}\cap V_{\Sigma_{1}} \cap
            \cdots \cap V_{\Sigma_{t}}
\end{equation}
has dimension zero: it is an oriented set of points, contained
entirely in the top stratum of the $V_{x_{i}}$. The next lemma is
proved by the same dimension-counting that is used for the case $N=2$:

\begin{lemma}
    Under the stated transversality assumptions, the zero-dimensional
    intersection \eqref{eq:q-intersection} is compact. \qed
\end{lemma}

Finally, we can state:

\begin{proposition}
     Fix $N$, and let $X$ be a closed, oriented smooth $4$-manifold with
    $b^{+}_{2}(X) \ge 2$, equipped with a homology orientation
    $o_{X}$.  Let $w$ be a line bundle with $c_{1}(w)$ coprime to
    $N$, and let $P$ be a $U(N)$ bundle on $X$ with
    $\det(P)=w$, such that the corresponding moduli space
    $M^{w}_{\kappa}$ has formal dimension an even integer $d$.
     Let $V_{x_{i}}$  and $V_{\Sigma_{j}}$ be as above. Then the signed
     count of the points in the intersection \eqref{eq:q-intersection}
     has the following properties:

     \begin{enumerate}
        \item It is independent of the choice of Riemannian metric $g$,
        perturbation $\bomega$, the neighborhoods $\nu(x_{i})$ and
        $\nu(\Sigma_{j})$, and the sections defining the submanifolds
        $V_{x_{i}}$ and $V_{\Sigma_{j}}$, subject to the conditions
        laid out in the construction above.

        \item It depends on the $\Sigma_{j}$ only through their
        homology classes $[\Sigma_{j}]$.

        \item It is linear in the homology class of each $\Sigma_{j}$.
     \end{enumerate}
\end{proposition}

\begin{proof}
    The proof can be modeled on the arguments given in
    \cite{Donaldson-polynomials} and \cite{Donaldson-Kronheimer}.
\end{proof}

Staying with our assumptions that $b_{2}^{+}(X)\ge 2$, that $c_{1}(w)$
is coprime to $N$ and that a homology orientation $o_{X}$ is given, we
now define the polynomial invariant
\[
                   q^{w}(X) : \bbA(X) \to \Q
\]
by declaring its value on the monomial $z$ given in
\eqref{eq:z-monomial} to be given by
\[
            q^{w}(X)(z) = (2N)^{-r-t} \# ( M^{w}_{\bomega,\kappa} \cap V_{x_{1}} \cap
            \cdots \cap V_{x_{r}}\cap V_{\Sigma_{1}} \cap
            \cdots \cap V_{\Sigma_{t}} ),
\]
where $\kappa$ is chosen (if possible) to make the moduli space
$d$-dimensional. If there is no such $\kappa$, we define $q^{w}(X)(z)$
to be zero for this monomial of degree $d$.
The factor $(2N)^{-r-t}$ is there because the subvarieties $V_{\Sigma}$ and
$V_{\mu}$ are dual to $(2N)\mu([\Sigma])$ and $(2N)\mu([x])$
respectively: so if the moduli space were actually compact, then the
$q^{w}(X)(z)$ would be the pairing
\[
         q^{w}(X)(z) = \bigl\langle \mu(z) , [M^{w}_{\bomega,\kappa}]
         \bigr\rangle.
\]

The definition of $q^{w}(X)$ is set up so that the integer invariant
that we defined earlier can be recovered from this more general
invariant as the value $q^{w}(X)(1)$.

\subsection{The non-coprime case}

Thus far, we have defined $q^{w}(X)$ only in the case that $c_{1}(w)$ is
coprime to $N$,  because this condition allows us to avoid the
difficulties arising from reducibles. In the case that $c_{1}(w)$ is
\emph{not} coprime to $N$, however, we can still give a useful
definition of $q^{w}(X)$ by a formal device previously used in
\cite{KM-Structure}.

Given $X$ and $w$ (with $b^{+}(X) \ge 2$ as
usual), we consider the blow-up
\[
            \tilde X = X \# \bar{\CP}^{2}.
\]
A maximal positive-definite subspace of $H^{2}(X;\R)$ is also a
maximal positive-definite subspace of $H^{2}(\tilde X;\R)$ (which is
canonically a direct sum); so a homology-orientation $o_{X}$
determines a homology-orientation $o_{\tilde X}$ for the blow-up. We
set
\[
                \tilde w = w + e,
\]
where $e = \mathrm{PD}(E)$ is the Poincar\'e dual of a chosen
generator of $E$ in $H_{2}(\bar{CP}^{2};\Z)$. Notice that $\tilde w$
is coprime to $N$ on $\tilde X$, so $q^{\tilde w}(\tilde X)$ is
already well-defined. We now \emph{define}
$q^{w}(X)$ in terms of $q^{\tilde w}(\tilde X)$ by the formula
\begin{equation}\label{eq:blow-up}
          q^{w}(X)(z) =   q^{\tilde w}(\tilde X)(E^{N-1}z), 
\end{equation}
for $z \in \bbA(X) \subset \bbA(\tilde X)$. (The result is independent
of the choice of generator $E$.)

The rationale for this definition, as in \cite{KM-Structure}, is that
the formula \eqref{eq:blow-up} is an \emph{identity} in the case that
$c_{1}(w)$ is already coprime to $N$ (so that the left-hand side is
defined by our previous construction). This is a blow-up
formula, which in the case $N=2$ was first proved in \cite{Kotschick}.
We now outline the proof of this blow-up formula, for general $N$.

On the manifold $\bar\CP^{2}$, let $Q$ be the $U(N)$ bundle with
$c_{1}(Q) = e$ and $c_{2}(Q) = 0$. The characteristic number
$\kappa(Q)$
for this bundle is $(N-1)/2N$, so the formal dimension $d$ of the
moduli space is $2(N-1) - (N^{2}-1)$. If we pick a base-point $y_{0}
\in \bar\CP^{2}$ and consider the \emph{framed} moduli space (the
quotient of the space of solutions by the based gauge group), then the
formal dimension is $2(N-1)$. The fact that $\bar\CP^{2}$ has positive
scalar curvature and anti-self-dual Weyl curvature means that the
framed moduli space is regular; and 
since its dimension is less than $N^{2}-1$, the
moduli space consists only of reducibles. The framed moduli space
therefore consists of a single $\CP^{N-1}$, which is the gauge orbit
of a connection $[A]$ compatible with the splitting of $Q$ as a sum 
\[
            Q = l \oplus Q'.
\]
(The first is the line bundle with $c_{1}(l) = e$, and the second
summand is a trivial $U(N-1)$ bundle.) The stabilizer $\Gamma_{A}$ in
the gauge group $\G$ is $S(U(1) \times U(N-1))$. Let us write $N$ for
this framed moduli space. There is a universal bundle
\[
           \ad( \mathbb{Q}) \to \bar\CP^{2}\times N
\]
and a corresponding slant-product map
\[
            \mu : H_{i}(\bar{\CP}^{2}) \to H^{4-i}(N).
\]
We can describe $\mathbb{Q}$ as the quotient
\[
           \ad( \mathbb{Q}) =
           \frac{\ad (Q) \times \SU(N) }
           { S(U(1) \times U(N-1))},
\]
which makes it the adjoint bundle of the $U(N)$ bundle
\[
         \mathbb{Q}=
           \frac{Q \times U(N) }
           { U(1) \times U(N-1)}.
\]
The decomposition of $Q$ gives rise to a decomposition
\[
            \mathbb{Q} = (w\boxtimes \tau) \oplus (\C \boxtimes
            \tau^{\bot}),
\]
where $\tau \to N = \CP^{N-1}$ is the tautological line bundle and
$\tau^{\bot}$ is its complement. On the submanifold $\CP^{1}\times
\CP^{1} \subset \bar{\CP}^{2}\times N$, the bundle $\mathbb{Q}$ has
$c_{1}^{2} = 0$ and $c_{2}$ is a generator. So $\mu(E)$ is an
integral
generator
for $H^{2}(N) = \Z$. In particular,
\begin{equation}\label{eq:N-mumap}
        \bigl\langle \mu(E)^{N-1}, [N] \bigr\rangle = \pm 1.
\end{equation}

Now we consider the blow-up formula in the case that $w$ is coprime to
$N$ on $X$ and that the moduli space $M^{w}_{\bomega, \kappa}(X)$ is
zero-dimensional, leading to an integer invariant $q^{w}(X)$.  Let
$x_{0}$ be a base-point in $X$. Using Lemma~\ref{lem:compactSupport}
and the conformal invariance of the anti-self-duality equations, we
can find a good pair $(g,\bomega)$ with the property that the support
of the perturbation is disjoint from a neighborhood of $x_{0}$. We can
then form a connected sum $\tilde X = X \# \bar{\CP}^{2}$, attaching
the $\bar\CP^{2}$ summand to $X$ in a neighborhood of $x_{0}$. We
equip the connected sum with a metric $g_{R}$ containing a neck of
length $R$, in the usual way. On $\tilde X$, we form the bundle
$\tilde P$ as the sum of the bundles $P\to X$ and $Q \to \bar\CP^{2}$.
As in \ref{subsec:gluing}, we may regard $\bomega$ as defining a
perturbation on $\tilde X$, with support disjoint from the neck and
from the $\bar\CP^{2}$ summand. The formal dimension of the moduli
space on $\tilde X$ is $2(N-1)$, and the moduli space is compact once
$R$ is sufficiently large. Standard gluing techniques apply, to show
that the moduli space for large $R$ is a product of the
zero-dimensional moduli space $M^{w}_{\bomega,\kappa}(X)$ and the
framed moduli space $N = \CP^{N-1}$ for $\bar\CP^{2}$.  This is a
diffeomorphism of oriented manifolds, up to an overall sign. From
\eqref{eq:N-mumap}, we deduce that
\[
            q^{w}(X)(1) = \pm q^{\tilde w}(\tilde X)(E^{N-1}),
\]
which is the desired formula \eqref{eq:blow-up}, up to sign, for the
case $z=1$. To prove the general case, we  take $z$ to be a
monomial as in \eqref{eq:z-monomial}, and we replace the moduli space
$M^{w}_{\bomega,\kappa}(X)$ by the cut-down moduli space
\eqref{eq:q-intersection}, having chosen the neighborhoods
$\nu(\Sigma_{j})$ etc.~to be disjoint from the region in which the
connected sum is formed. The argument is then essentially the same as
the case $z=1$. To identify the correct sign for the blow-up formula,
we can again make a comparison with the K\"ahler case.

Having dealt in this way with the non-coprime case,
we are allowed now to take $w$ trivial, so as to arrive at a
satisfactory definition of the $\SU(N)$ Donaldson invariant
corresponding to a bundle $P$ with $c_{1}=0$.

\section{Further thoughts}

\subsection{Higher Pontryagin classes}

We have defined $\mu$ in \eqref{eq:mu-map} using the $4$-dimensional
characteristic class $\cc$, taking our lead from the usual definition
of the $N=2$ invariants. But there are other possibilities. The
easiest generalizations is of the $4$-dimensional class $\mu([x])$,
where $[x]$ is the generator of $H_{0}(X)$. Recall that the cohomology
class $\mu([x])$ is a rational multiple of the first Pontryagin class,
of the base-point bundle $\ad(\mathbb{P}_{x})\to \B^{*}$. Let us drop
the factor $-(1/2N)$ in from the definition and write
\[
                \nu_{4} =  p_{1}(\ad(\mathbb{P}_{x})).
\]
We can consider more generally now the class
\[
  \nu_{4j} =  p_{j}(\ad(\mathbb{P}_{x}))
\]
in $H^{4j}(\B^{*})$. On the moduli spaces $M^{w}_{\bomega,\kappa}$,
this can be represented as the dual of a subvariety
\[
                M^{w}_{\bomega,\kappa} \cap V_{j,x},
\]
where $V_{j,x}$ is again pulled back from a neighborhood of $x$. We
restrict as before to the coprime case, with $b^{+}_{2}(X) \ge 2$, in
which case we consider a zero-dimensional intersection
\begin{equation}\label{eq:higher-pj}
                M^{w}_{\bomega,\kappa} \cap V_{j_{1},x_{1}}
                \cap \dots \cap V_{j_{r},x_{r}}.
\end{equation}
If $j_{i} \le N-1$ for all $i$, and if the $V_{j_{i},x_{i}}$ are
chosen so that all relevant intersections are transverse, then
this zero-dimensional intersection will be compact.  So we can define
Donaldson-type invariants using the characteristic class $\nu_{j}$ of
the base-point fibration, as long as $j \le N-1$.
 
For $j \ge N$, the dimension-counting fails. For example, suppose that
$M^{w}_{\bomega,\kappa}$ is $4N$-dimensional, so that
$M^{w}_{\bomega,\kappa-1}$ is $0$-dimensional, and suppose that  $w$
is coprime to $N$ and that the moduli spaces are regular and without
reducibles. We can try taking $j = N$, and considering the
$0$-dimensional intersection
\[
            M^{w}_{\bomega,\kappa} \cap V_{N, x},
\]
so as to arrive at an evaluation of $p_{N}$ of the base-point
fibration on the (potentially non-compact) moduli space. The problem
is that the above zero-dimensional intersection need not be compact.
We cannot rule out the possibility that there is a sequence $[A_{n}]$
of solutions belonging to this intersection, converging to an ideal
connection $([A'], x')$, with $[A'] \in M^{w}_{\bomega,\kappa-1}$ and
$x'$ in the neighborhood $\nu(x)$.

\subsection{Other Lie groups}

The approach we have taken to define the $\SU(N)$ polynomial
invariants (even in the case of a bundle with $c_{1}=0$) makes
essential use of a $U(N)$ bundle on the blow-up, in order to avoid
reducible solutions. This approach does not seem to be available in
general: we cannot define $E_{8}$ Donaldson invariants this way, for
example. It is possible that by taking the instanton number to be
large, one could sidestep the difficulties that arise from reducible
solutions, and so arrive at an alternative construction. This was
Donaldson's original approach, in \cite{Donaldson-polynomials}. It
seems likely that this approach, if successful, would be technically
more difficult than what we have done here for $\SU(N)$.

\subsection{Mahler measure}

The large-$N$ behaviour of the $\SU(N)$ Donaldson invariants is
briefly considered in \cite{Marino-Moore}. In the particular examples
we have studied in this paper, the leading term of the invariant for
large $N$ is easily visible, and turns out to have a rather
interesting interpretation.

We turn again to the manifold $Z_{K}$ with torus boundary, obtained by
taking the product of the circle and a knot-complement, and we look
again at the integer-valued invariant $q^{w}(Z_{K};\alpha_{0})$. The
formula for this invariant in
Proposition~\ref{prop:relative-calcuation} relates it to the geometric
mean of the Alexander polynomial $\Delta_{K}$ on the unit circle: as
long as there are no zeros of $\Delta_{K}$ on the unit circle, we have
\[
           q^{w}(Z_{K};\alpha_{0})_{N} 
           \sim \alpha^{N}, 
\]
as $N$ increases through odd integers, where
\[
                \alpha = \exp 
                \int_{0}^{1}\log | \Delta_{K}(e^{2\pi i t}) | \, dt.
\]
This asymptotic result is very easily obtained by considering 
the expression for
\[
            (1/N) \log | q^{w}(Z_{K};\alpha_{0})_{N} |
\]
arising from Proposition~\ref{prop:relative-calcuation}: it is a Riemann-sum
approximation to the integral of the smooth function $| \Delta_{K}(z)|$ on
the circle. The remaining point is then that the Riemann sum approximates
the integral of a periodic function with an error which is at most $O(1/N^{3})$ when the
integrand is $C^{2}$.

If $\Delta_{K}$ has zeros
on the unit circle, the situation is much more delicate. We have only the weaker statement
\begin{equation}\label{eq:delicate-Mahler}
      \lim_{N\to\infty}   (1/N)  \log  q^{w}(Z_{K};\alpha_{0})_{N} 
           = \log \alpha,
\end{equation}
as $N$ runs through those odd integers for which
$q^{w}(Z_{K};\alpha_{0})_{N}$ is non-zero.
The proof of \eqref{eq:delicate-Mahler} is given in \cite{GA-Short}
and in \cite{Riley}, in answer to a question raised by Gordon in
\cite{Gordon}. A different proof and an extension of the result are
contained in \cite{Silver-Williams-a}.

The quantity $\alpha$ is known as the \emph{Mahler measure} of
the polynomial. It can also be expressed as
\[
            \alpha = |a| \prod_{\lambda} \max (1, |\lambda|),
\]
where $a$ is the leading coefficient of the polynomial and $\lambda$
runs through the roots in the complex plane. The Mahler measure of the Alexander
polynomial, as well as the Mahler measures of the Jones and $A$-polynomials, have appeared
elsewhere in the literature, in connection with interesting quantities
such as the hyperbolic volume of the knot complement. See
\cite{Silver-Williams, Murakami} and \cite{Boyd}, for example.

\bibliography{higherrank}

\end{document}